\documentclass[12pt,a4paper,leqno]{article}

\usepackage{latexsym}
\usepackage{amssymb,exscale}
\usepackage[centertags]{amsmath}
\usepackage{amsthm}
\usepackage[all]{xy}
\usepackage{graphicx}

\usepackage[dvips]{hyperref}

\numberwithin{equation}{section}

\swapnumbers
\theoremstyle{definition}
\newtheorem{theorem}[equation]{Theorem}

\newtheorem{lemma}[equation]{Lemma}
\newtheorem{proposition}[equation]{Proposition}
\newtheorem{corollary}[equation]{Corollary}
\newtheorem{claim}[equation]{Claim}
\newtheorem{definition}[equation]{Definition}
\newtheorem{remark}[equation]{Remark}

\newtheorem{example}[equation]{Example}
\newtheorem{counterexample}[equation]{Counterexample}

\renewcommand{\phi}{\varphi}

\newcommand{\D}{\mathrm{d}}
\newcommand{\E}{\mathrm{e}}

\newcommand{\ti}{\tilde}

\renewcommand{\(}{\bigl(}
\renewcommand{\)}{\bigr)\vphantom{)}}

\newcommand{\equi}{\;\Longleftrightarrow\;}

\newcommand{\imply}{\;\;\,\Longrightarrow\;\;\,}

\newcommand{\imp}{$ \Longrightarrow $ }

\newcommand{\Poisson}{\operatorname{Poisson}}
\newcommand{\CS}{\operatorname{CS}}
\newcommand{\FCS}{\operatorname{FCS}}
\newcommand{\DCS}{\operatorname{DCS}}
\newcommand{\Exp}{\operatorname{Exp}}
\newcommand{\mes}{\operatorname{mes}}

\newcommand{\modO}{{\operatorname{mod}\,0}}

\newcommand{\OI}{\mbox{0\hspace{1pt}--1}}

\newcommand{\One}{\mathbf1}

\newcommand{\BB}{\mathsf B}

\newcommand{\eps}{\varepsilon}
\newcommand{\si}{\sigma}
\newcommand{\Si}{\Sigma}

\newcommand{\om}{\omega}
\newcommand{\Om}{\Omega}

\newcommand{\al}{\alpha}
\newcommand{\be}{\beta}
\newcommand{\X}{\mathcal X}
\newcommand{\Ec}{\mathcal E}
\newcommand{\F}{\mathcal F}

\newcommand{\A}{\mathcal A}
\newcommand{\B}{\mathcal B}

\newcommand{\const}{{\mathrm{const}}}
\newcommand{\Ex}{\mathbb E\,}
\renewcommand{\Pr}[1]{\mathbb{P}\mskip1.5mu\(\mskip1.5mu#1\mskip1.5mu\)}
\newcommand{\R}{\mathbb R}
\newcommand{\Q}{\mathbb Q}

\newcommand{\cE}[2]{\mathbb{E}\mskip1.5mu\(\mskip1.5mu#1\mskip1.5mu
 \big|\mskip1.5mu#2\mskip1.5mu\)}

\newcommand{\cP}[2]{\mathbb{P}\mskip1.5mu\(\mskip1.5mu#1\mskip1.5mu
 \big|\mskip1.5mu#2\mskip1.5mu\)}

\newcommand{\sif}{$\sigma$\nobreakdash-field}

\newcommand{\finite}[1]{$#1$\nobreakdash-\hspace{0pt}finite}
\newcommand{\based}[1]{$#1$\nobreakdash-\hspace{0pt}based}

\newcommand{\almost}[1]{$#1$\nobreakdash-\hspace{0pt}almost}
\newcommand{\invariant}[1]{$#1$\nobreakdash-\hspace{0pt}invariant}
\newcommand{\valued}[1]{$#1$\nobreakdash-\hspace{0pt}valued}

\newcommand{\measurable}[1]{$#1$\nobreakdash-\hspace{0pt}measurable}

\begin{document}

\title{Brownian local minima,\\
 random dense countable sets\\
 and random equivalence classes} 

\author{Boris Tsirelson}

\date{}
\maketitle

\stepcounter{footnote}
\footnotetext{%
 This research was supported by \textsc{the israel science foundation}
 (grant No.~683/05).}

\begin{abstract}
A random dense countable set is characterized (in distribution) by
independence and stationarity. Two examples are \emph{Brownian local
minima} and \emph{unordered infinite sample}. They are identically
distributed. A framework for such concepts, proposed here, includes a
wide class of random equivalence classes.
\end{abstract}

\section*{Introduction}
Random countable sets arise naturally from the Brownian motion (local
extrema, see \cite[2.9.12]{KS}), percolation (double, or four-arm points,
see \cite{CFN}), oriented percolation (points of type $(2,1)$, see
\cite[Th.~5.15]{FIN}) etc. They are scarcely investigated, because they
fail to fit into the usual framework. They cannot be treated as random
elements of `good' (Polish, standard) spaces. The framework of adapted
Poisson processes, used by Aldous and Barlow \cite{AB}, does not apply
to the Brownian motion, since the latter cannot be correlated with a
Poisson process adapted to the same filtration. The `hit-and-miss'
framework used by Kingman \cite{Ki} and Kendall \cite{Ke00} fails to
discern the clear-cut distinction between Brownian local minima and,
say, randomly shifted set of rational numbers. A new approach introduced
here catches this distinction, does not use adapted processes, and shows
that Brownian local minima are distributed like an infinite sample in
the following sense (see Theorem \ref{6n10}).

{
\newtheorem*{thrm}{Theorem}
\begin{thrm}
\begin{sloppypar}
There exists a probability measure $ P $ on the product space $ C[0,1]
\times (0,1)^\infty $ such that 
\end{sloppypar}

(a) the first marginal of $ P $ (that is, projection to the first
factor) is the Wiener measure on the space $ C[0,1] $ of continuous
paths $ w : [0,1] \to \R $;

(b) the second marginal of $ P $ is the Lebesgue measure on the
cube $ (0,1)^\infty $ of infinite (countable) dimension;

(c) \almost{P} all pairs $ (w,u) $, $ w \in C[0,1] $, $ u =
(u_1,u_2,\dots) \in (0,1)^\infty $, are such that the numbers $
u_1,u_2,\dots $ are an enumeration of the set of all local minimizers of
the Brownian path $ w $.
\end{thrm}
}

Thus, the conditional distribution of $ u_1,u_2,\dots $ given $ w $
provides a (randomized) enumeration of Brownian minimizers by
independent uniform random variables.

The same result holds for every random dense countable set that
satisfies conditions of independence and stationarity, see Definitions
\ref{4n2}, \ref{6n7} and Theorem \ref{6n8}. Two-dimensional
generalizations, covering the percolation-related models, are possible.

On a more abstract level the new approach is formalized in Sections
\ref{sect7}, \ref{sect8} in the form of `borelogy' that combines some
ideas of descriptive set theory \cite{Ke99} and diffeology
\cite{IZ}. Random elements of various quotient spaces fit into the new
framework. Readers that like abstract concepts may start with these
sections.

\section[]{\raggedright Main lemma}
\label{sect1}Before the main lemma we consider an instructive special case.

\begin{example}
Let $ \mu $ be the uniform distribution on the interval $ (0,1) $ and $
\nu $ the triangle distribution on the same interval, that is,
\[
\mu(B) = \int_B \D x \, , \qquad \nu(B) = \int_B 2x \, \D x
\]
for all Borel sets $ B \subset (0,1) $. On the space $ (0,1)^\infty $ of
sequences, the product measure $ \mu^\infty $ is the joint
distribution of uniform i.i.d.\ random variables, while $ \nu^\infty $
is the joint distribution of triangular i.i.d.\ random variables. I
claim existence of a probability measure $ P $ on $ (0,1)^\infty \times
(0,1)^\infty $ such that

(a) the first marginal of $ P $ is equal to $ \mu^\infty $ (that is, $ P
( B \times (0,1)^\infty ) = \mu^\infty (B) $ for all Borel sets $ B
\subset (0,1)^\infty $);

(b) the second marginal of $ P $ is equal to $ \nu^\infty $ (that is, $ P
( (0,1)^\infty \times B ) = \nu^\infty (B) $ for all Borel sets $ B
\subset (0,1)^\infty $);

(c) \almost{P} all pairs $ (x,y) $, $ x=(x_1,x_2,\dots) \in (0,1)^\infty
$, $ y=(y_1,y_2,\dots) \in (0,1)^\infty $ are such that
\[
\{ x_1, x_2, \dots \} = \{ y_1, y_2, \dots \} \, ;
\]
in other words, the sequence $ y $ is a permutation of the sequence $ x
$. (A random permutation, of course.)

A paradox: the numbers $ y_k $ are biased toward $ 1 $, the numbers $
x_k $ are not; nevertheless they are just the same numbers!

\begin{figure}
\setlength{\unitlength}{1cm}
\begin{picture}(15,5.9)
\put(0.0,0.8){\includegraphics{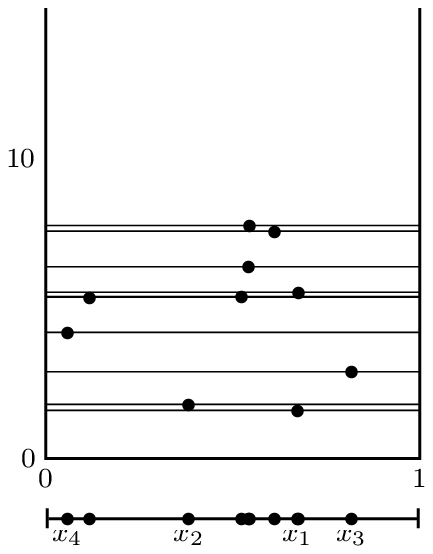}}
\put(4.7,0.61){\includegraphics{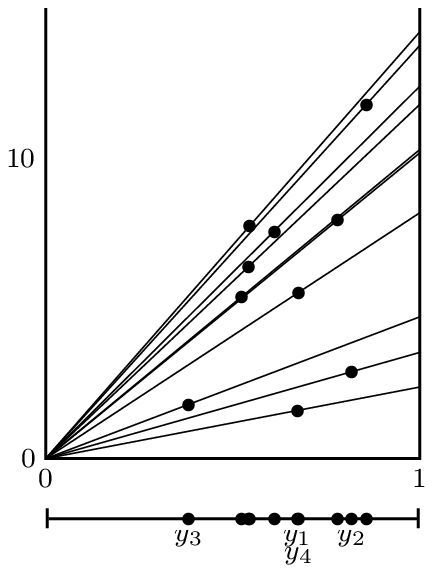}}
\put(9.4,0.63){\includegraphics{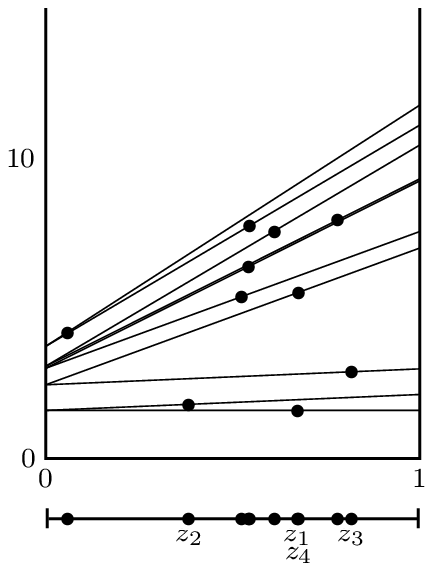}}
\put(2,0.05){(a)}
\put(6.6,0.05){(b)}
\put(11.2,0.05){(c)}
\end{picture}
\caption[]{\label{fig1}\small
different enumerations turn the same set into an infinite sample from:
(a) the uniform distribution, (b) the triangular distribution, (c) their
mixture (uniform $ z_1,z_3,\dots $ but triangular $ z_2,z_4,\dots $).}
\end{figure}

An explanation (and a sketchy proof) is shown on Fig.~\ref{fig1}(a,b). A
countable subset of the strip $ (0,1) \times (0,\infty) $ is a
realization of a Poisson point process. (The mean number of points in
any domain is equal to its area.) The first $ 10 $ points of the same
countable set are shown on both figures, but on Fig.~\ref{fig1}(a) the
points are ordered according to the vertical coordinate, while on
Fig.~\ref{fig1}(b) they are ordered according to the ratio of the two
coordinates. We observe that $ \{y_1,\dots,y_{10}\} $ is indeed biased
toward $ 1 $, while $ \{x_1,\dots,x_{10}\} $ is not. On the other hand,
$ y $ is a permutation of $ x $. (This time, $ y_1=x_1 $, $ y_2=x_3 $, $
y_3=x_2 $, \dots)
\end{example}

A bit more complicated ordering, shown on Fig.~\ref{fig1}(c), serves the
measure $ \mu \times \nu \times \mu \times \nu \times \dots $, the joint
distribution of independent, differently distributed random variables.

In every case we use an increasing sequence of (random) functions $ h_n
: (0,1) \to [0,\infty) $ such that for each $ n $ the graph of $ h_n $
contains a Poisson point, while the region between the graphs of $
h_{n-1} $ and $ h_n $ does not. The differences $ h_n - h_{n-1} $ are
constant functions on Fig.~\ref{fig1}(a), triangular (that is, $ x
\mapsto \const \cdot x $) on Fig.~\ref{fig1}(b), while on
Fig.~\ref{fig1}(c) they are constant for odd $ n $ and triangular for
even $ n $.

Moreover, the same idea works for dependent random variables. In this
case $ h_n - h_{n-1} $ is proportional to the conditional density, given
the previous points. We only need existence of conditional densities and
divergence of their sum (in order to exhaust the strip).

Here is the main lemma.

\begin{lemma}\label{3.2}
Let $ \mu $ be a probability measure on $ (0,1)^\infty $ such that

(a) for every $ n $ the marginal distribution of the first $ n $
coordinates is absolutely continuous;

(b) for almost all $ x \in (0,1) $ and \almost{\mu} all $
(x_1,x_2,\dots) \in (0,1)^\infty $,
\[
\sum_{n=1}^\infty \frac{ f_{n+1}(x_1,\dots,x_n,x) }{ f_n(x_1,\dots,x_n)
} = \infty \, ;
\]
here $ f_n $ is the density of the first $ n $ coordinates.

Let $ \nu $ be another probability measure on $ (0,1)^\infty $
satisfying the same conditions (a), (b). Then there exists a probability
measure $ P $ on $ (0,1)^\infty \times (0,1)^\infty $ such that

(c) the first marginal of $ P $ is equal to $ \mu $;

(d) the second marginal of $ P $ is equal to $ \nu $;

(e) \almost{P} all pairs $ (x,y) $, $ x=(x_1,x_2,\dots) \in (0,1)^\infty
$, $ y=(y_1,y_2,\dots) \in (0,1)^\infty $ are such that
\[
\{ x_1, x_2, \dots \} = \{ y_1, y_2, \dots \} \, .
\]
\end{lemma}

In other words, the sequence $ y $ is a permutation of the sequence $ x
$ (since $ x_k $ are pairwise different due to absolute continuity, as
well as $ y_k $).

The rest of the section is occupied by the proof of the main lemma.

Throughout the proof, Poisson point processes on the strip $ (0,1)
\times [0,\infty) $ (or its measurable part) are such that the mean
number of points in any measurable subset is equal to its
two-dimensional Lebesgue measure. Random variables (and random
functions) are treated here as measurable functions of the original
Poisson point process on the strip.

We start with three rather general claims.

\begin{claim}\label{3.3}
(a) A Poisson point process on the strip $ (0,1) \times [0,\infty) $
may be treated as the set $ \Pi $ of random points $ (U_n,T_1+\dots+T_n)
$ for $ n=1,2,\dots $, where $ U_1,T_1,U_2,T_2,\dots $ are independent
random variables, each $ U_k $ is distributed uniformly on $ (0,1) $,
and each $ T_k $ is distributed $ \Exp(1) $ (that is, $ \Pr{ T_k > t } =
\E^{-t} $ for $ t \ge 0 $);

(b) conditionally on $ (U_1,T_1) $, the set $ \Pi_1 = \{
(U_n,T_1+\dots+T_n) : n \ge 2 \} $ is (distributed as) a Poisson point
process on $ (0,1) \times [T_1,\infty) $.
\end{claim}

The proof is left to the reader.

\begin{claim}\label{3.4}
Let $ f : (0,1) \to [0,\infty) $ be a measurable function satisfying $
\int_0^1 f(x) \, \D x = 1 $, and $ \Pi $ be a Poisson point process on $
(0,1) \times [0,\infty) $. Then

(a) the minimum
\[
t_1 = \min_{(x,y)\in\Pi} \frac{y}{f(x)}
\]
(where $ y/0=\infty $) is reached at a single point $ (x_1,y_1) \in \Pi
$;

(b) $ x_1 $ and $ t_1 = y_1 / f(x_1) $ are independent, $ t_1 $ is
distributed $ \Exp(1) $, and the distribution of $ x_1 $ has the density
$ f $;

(c) conditionally on $ (x_1,y_1) $, the set $ \Pi_1 = \Pi \setminus
\{(x_1,y_1)\} $ is (distributed as) a Poisson point process on $ \{
(x,y) : 0<x<1, \, t_1 f(x) < y < \infty \} $.
\end{claim}

\begin{proof}
A special case, $ f(x)=1 $ for all $ x $, follows from \ref{3.3}.

A more general case, $ f(x)>0 $ for all $ x $, results from the special
case by the transformation $ (x,y) \mapsto (F(x), y/f(x)) $ where $ F(x)
= \int_0^x f(x') \, \D x' $. The transformation preserves Lebesgue
measure on $ (0,1) \times [0,\infty) $, therefore it preserves also the
Poisson point process.

In the general case the same transformation sends $ A \times [0,\infty)
$ to $ (0,1) \times [0,\infty) $, where $ A = \{ x : f(x)>0 \} $.
Conditionally on $ (x_1,t_1) $ we get a Poisson point process on $ \{ 
(x,y) : x \in A, \, t_1 f(x)<y<\infty \} $ independent of the Poisson
point process on $ \{ (x,y) : x \notin A, \, 0<y<\infty \} $.
\end{proof}

\begin{claim}\label{3.5}
Let $ f : (0,1) \to [0,\infty) $ and $ g : (0,1) \times (0,1) \to
[0,\infty) $ be measurable functions satisfying $ \int_0^1 f(x) \, \D x
= 1 $ and $ \int_0^1 g(x_1,x_2) \, \D x_2 = 1 $ for almost all $ x_1
$. Let $ \Pi $ be a Poisson point process on $ (0,1) \times [0,\infty)
$, while $ (x_1,y_1) $, $ t_1 $ and $ \Pi_1 $ be as in \ref{3.4}. Then

(a) the minimum
\[
t_2 = \min_{(x,y)\in\Pi_1} \frac{y-t_1 f(x)}{g(x_1,x)}
\]
is reached at a single point $ (x_2,y_2) \in \Pi_1 $;

(b) conditionally on $ (x_1,y_1) $ we have: $ x_2 $ and $ t_2 = (y_2 -
t_1 f(x_2) ) / g(x_1,x_2) $ are independent, $ t_2 $ is distributed $
\Exp(1) $, and the distribution of $ x_2 $ has the density $
g(x_1,\cdot) $ (conditional distributions are meant);

(c) conditionally on $ (x_1,y_1), (x_2,y_2) $, the set $ \Pi_2 = \Pi
\setminus \{(x_1,t_1), (x_2,y_2)\} $ is a Poisson point process on $ \{
(x,y) : 0<x<1, t_1 f(x) + t_2 g(x_1,x) < y < \infty \} $. 
\end{claim}

\begin{proof}
After conditioning on $ (x_1,y_1) $ we apply \ref{3.4} to the Poisson
point process $ \{ (x,y-t_1 f(x)) : (x,y) \in \Pi_1 \} $ on $ (0,1)
\times [0,\infty) $ and the function $ g(x_1,\cdot) : (0,1) \to
[0,\infty) $.
\end{proof}

Equipped with these claims we prove the main lemma as
follows. Introducing conditional densities
\[
g_{n+1} ( x | x_1,\dots,x_n ) = \frac{ f_{n+1} (x_1,\dots,x_n,x) }{
f_n(x_1,\dots,x_n) } \, , \qquad g_1(x) = f_1(x)
\]
and a Poisson point process $ \Pi $ on the strip $ (0,1) \times
[0,\infty) $, we construct a sequence of points $ (X_n,Y_n) $ of $ \Pi
$, random variables $ T_n $ and random functions $ H_n $ as follows:
\begin{align*}
T_1 &= \min_{(x,y)\in\Pi} \frac{y}{g_1(x)} = \frac{Y_1}{g_1(X_1)} \, ; \\
H_1(x) &= T_1 g_1(x) \, ; \\
T_2 &= \min_{(x,y)\in\Pi_1} \frac{y-H_1(x)}{g_2(x|X_1)} =
 \frac{Y_2-H_1(X_2)}{g_2(X_2|X_1)} \, ; \\
H_2(x) &= H_1(x) + T_2 g_2(x|X_1) \, ; \\
& \quad \dots \\
T_n &= \min_{(x,y)\in\Pi_{n-1}}
 \frac{y-H_{n-1}(x)}{g_n(x|X_1,\dots,X_{n-1})} =
 \frac{Y_n-H_{n-1}(X_n)}{g_n(X_n|X_1,\dots,X_{n-1})} \, ; \\
H_n(x) &= H_{n-1}(x) + T_n g_n(x|X_1,\dots,X_{n-1}) \, ; \\
& \quad \dots
\end{align*}
here $ \Pi_n $ stands for $ \Pi \setminus \{ (X_1,Y_1),\dots,(X_n,Y_n)
\} $. By \ref{3.4}(b), $ X_1 $ and $ T_1 $ are independent, $ T_1 \sim
\Exp(1) $ and $ X_1 \sim g_1 $. By \ref{3.4}(c) and \ref{3.5}(b),
conditionally on $ X_1 $ and $ T_1 $, $ \Pi_1 $ is a Poisson point
process on $ \{ (x,y) : y > H_1(x) \} $, while $ X_2 $ and $ T_2 $ are
independent, $ T_2 \sim \Exp(1) $ and $ X_2 \sim g_2(\cdot|X_1) $. It
follows that $ T_1 $, $ T_2 $ and $ (X_1,X_2) $ are independent, and the
joint distribution of $ X_1,X_2 $ has the density $ g_1(x_1)
g_2(x_2|x_1) = f_2(x_1,x_2) $. By \ref{3.5}(c), $ \Pi_2 $ is a Poisson
point process on $ \{ (x,y) : y > H_2(x) \} $ conditionally, given $
(X_1,Y_1) $ and $ (X_2,Y_2) $.

The same arguments apply for any $ n $. We get two independent sequences,
$ (T_1,T_2,\dots) $ and $ (X_1,X_2,\dots) $. Random variables $ T_n $
are independent, distributed $ \Exp(1) $ each. The joint distribution of
$ X_1, X_2, \dots $ is equal to $ \mu $, since for every $ n $ the joint
distribution of $ X_1,\dots,X_n $ has the density $ f_n $. Also, $ \Pi_n
$ is a Poisson point process on $ \{ (x,y) : y > H_n(x) \} $
conditionally, given $ (X_1,Y_1), \dots, (X_n,Y_n) $.

\begin{claim}\label{3.6}
$ H_n(x) \uparrow \infty $ for almost all pairs $ (\Pi,x) $.
\end{claim}

\begin{proof}
By \ref{3.2}(b), $ \sum_n g_{n+1} (x|x_1,\dots,x_n) = \infty $ for
almost all $ x \in (0,1) $ and \almost{\mu} all $ (x_1,x_2,\dots)
$. Therefore $ \sum_n g_{n+1} (x|X_1,\dots,X_n) = \infty $ for
almost all $ x \in (0,1) $ and almost all $ \Pi $. It is easy to see
that $ \sum_n c_n T_n = \infty $ a.s.\ for each sequence $ (c_n)_n $
such that $ \sum_n c_n = \infty $. Taking into account that $
(T_1,T_2,\dots) $ is independent of $ (X_1,X_2,\dots) $ we conclude that
$ \sum_n T_n g_n (x|X_1,\dots,X_{n-1}) = \infty $ for almost all $ x $
and $ \Pi $. The partial sums of this series are nothing but $ H_n(x) $.
\end{proof}

Still, we have to prove that the points $ (X_n,Y_n) $ exhaust the set $
\Pi $. Of course, a non-random negligible set does not intersect $ \Pi $
a.s.; however, the negligible set $ \{ x : \lim_n H_n(x) < \infty \} $
is random.

\begin{claim}
The set $ \cap_n \Pi_n $ is empty a.s.
\end{claim}

\begin{proof}
It is sufficient to prove that $ \cap_n \Pi_{n,M} = \emptyset $ a.s.\
for every $ M \in (0,\infty) $, where $ \Pi_{n,M} = \{ (x,y) \in \Pi_n :
y < M \} $. Conditionally, given $ (X_1,Y_1), \dots, \linebreak[0]
(X_n,Y_n) $, the set $ \Pi_{n,M} $ is a Poisson point process on $ \{
(x,y) : H_n(x) < y < M \} $; the number $ |\Pi_{n,M}| $ of points in $
\Pi_{n,M} $ satisfies
\[
\cE{ |\Pi_{n,M}| }{ (X_1,Y_1), \dots, (X_n,Y_n) } = \int_0^1
(M-H_n(x))^+ \, \D x \, .
\]
Therefore
\[
\Ex |\Pi_{n,M}| = \Ex \int_0^1 (M-H_n(x))^+ \, \D x \, .
\]
By \ref{3.6} and the monotone convergence theorem, $ \Ex \int_0^1
(M-H_n(x))^+ \, \D x \to 0 $ as $ n \to \infty $. Thus, $ \lim_n
|H_{n,M}| = 0 $ a.s.
\end{proof}

\begin{sloppypar}
Now we are in position to finish the proof of the main lemma. Applying
our construction twice (for $ \mu $ and for $ \nu $) we get two
enumerations of a single Poisson point process on the strip,
\[
\{ (X_n,Y_n) : n=1,2,\dots \} = \Pi = \{ (X'_n,Y'_n) : n=1,2,\dots \} \,
,
\]
such that the sequence $ (X_1,X_2,\dots) $ is distributed $ \mu $ and
the sequence $ (X'_1,X'_2,\dots) $ is distributed $ \nu $. The joint
distribution $ P $ of these two sequences satisfies the conditions
\ref{3.2}(c,d,e).
\end{sloppypar}

\section[]{\raggedright Random countable sets}
\label{sect2}Following the `constructive countability' approach of Kendall
\cite[Def.~3.3]{Ke00} we treat a random countable subset of $ (0,1) $ as
\begin{equation}\label{2n0}
\om \mapsto \{ X_1(\om), X_2(\om), \dots \}
\end{equation}
where $ X_1, X_2, \dots : \Om \to (0,1) $ are random variables. (To be
exact, it would be called a random finite or countable set, since $
X_n(\om) $ need not be pairwise distinct.) It may happen that $ \{
X_1(\om), X_2(\om), \dots \} = \{ Y_1(\om), Y_2(\om), \dots \} $ for
almost all $ \om $ (the sets are equal, multiplicity does not matter);
then we say that the two sequences $ (X_k)_k $, $ (Y_k)_k $ of random
variables represent the same random countable set, and write $ \{
X_1,X_2,\dots \} = \{ Y_1,Y_2,\dots \} $. On the other hand it may
happen that the joint distribution of $ X_1,X_2,\dots $ is equal to the
joint distribution of some random variables $ X'_1,X'_2,\dots : \Om' \to
(0,1) $ on some probability space $ \Om' $; then we may say that the two
random countable sets $ \{ X_1,X_2,\dots \} $, $ \{ X'_1,X'_2,\dots \} $
are identically distributed. We combine these two ideas as follows. 

\begin{definition}\label{2m1}
\begin{sloppypar}
Two random countable sets $ \{ X_1,X_2,\dots \} $, $ \{ Y_1,Y_2,\dots
\} $ are \emph{identically distributed} (in other words, $ \{ Y_1,Y_2,\dots
\} $ is distributed like $ \{ X_1,X_2,\dots \} $), if there exists a
probability measure $ P $ on the space $ (0,1)^\infty \times
(0,1)^\infty $ such that
\end{sloppypar}

(a) the first marginal of $ P $ is equal to the joint distribution of $
X_1,X_2,\dots $;

(b) the second marginal of $ P $ is equal to the joint distribution of $
Y_1,Y_2,\dots $;

(c) \almost{P} all pairs $ (x,y) $, $ x=(x_1,x_2,\dots) \in (0,1)^\infty
$, $ y=(y_1,y_2,\dots) \in (0,1)^\infty $ are such that
\[
\{ x_1, x_2, \dots \} = \{ y_1, y_2, \dots \} \, .
\]
\end{definition}

A sufficient condition is given by Main lemma \ref{3.2}: if Conditions
(a), (b) of Main lemma are satisfied both by the joint distribution of $
X_1,X_2,\dots $ and by the joint distribution of $ Y_1,Y_2,\dots $ then
$ \{ X_1,X_2,\dots \} $ and $ \{ Y_1,Y_2,\dots \} $ are identically
distributed.

\begin{remark}
The relation defined by \ref{2m1} is transitive. Having a joint
distribution of two sequences $ (X_k)_k $ and $ (Y_k)_k $ and a joint
distribution of $ (Y_k)_k $ and $ (Z_k)_k $ we may construct an
appropriate joint distribution of three sequences $ (X_k)_k $, $ (Y_k)_k
$ and $ (Z_k)_k $; for example, $ (X_k)_k $ and $ (Z_k)_k $ may be made
conditionally independent given $ (Y_k)_k $.
\end{remark}

\begin{definition}\label{2n2}
\begin{sloppypar}
A random countable set $ \{ X_1,X_2,\dots \} $ has \emph{the uniform
distribution} (in other words, \emph{is uniform}), if $ \{ X_1,X_2,\dots
\} $ and $ \{ Y_1,Y_2,\dots \} $ are identically distributed for some
(therefore, all) $ Y_1,Y_2,\dots $ whose joint distribution satisfies
Conditions (a), (b) of Main lemma \ref{3.2}.
\end{sloppypar}
\end{definition}

The joint distribution of $ X_1, X_2, \dots $ may violate Condition
\ref{3.2}(b), see \ref{9n75}.

If $ X_1,X_2,\dots $ are i.i.d.\ random variables then the random
countable set $ \{ X_1,X_2,\dots \} $ may be called an \emph{unordered
infinite sample} from the corresponding distribution. If the latter has a
non-vanishing density on $ (0,1) $ then $ \{ X_1,X_2,\dots \} $ has the
uniform distribution. A paradox: the distribution of the sample does not
depend on the underlying one-dimensional distribution! (See also
\ref{9n8}.)

\begin{remark}\label{2m5}
If the joint distribution of $ X_1,X_2,\dots $ satisfies Condition (a)
of Main lemma (but may violate (b)) then the random countable set $ \{
X_1,X_2,\dots \} $ may be treated as a part (subset) of a uniform random
countable set, in the following sense.

We say that $ \{ X_1,X_2,\dots \} $ is distributed as a part of $ \{
Y_1,Y_2,\dots \} $ if there exists $ P $ satisfying (a), (b) of
\ref{2m1} and (c) modified by replacing the equality $ \{ x_1, x_2,
\dots \} = \{ y_1, y_2, \dots \} $ with the inclusion $ \{ x_1, x_2,
\dots \} \subset \{ y_1, y_2, \dots \} $.

Indeed, the proof of Main lemma uses Condition (b) only for exhausting
all elements of the Poisson random set.
\end{remark}

\section[]{\raggedright Selectors}
\label{sect3}A single-element part of a random countable set is of special interest.

\begin{definition}\label{3n1}
A \emph{selector} of a random countable set $ \{ X_1,X_2,\dots \} $ is a
probability measure $ P $ on the space $ (0,1)^\infty \times (0,1) $
such that

(a) the first marginal of $ P $ is equal to the joint distribution of $
X_1,X_2,\dots $;

(b) \almost{P} all pairs $ (x,z) $, $ x=(x_1,x_2,\dots) \in (0,1)^\infty
$, $ z \in (0,1) $ satisfy
\begin{equation}\label{3n2}
z \in \{ x_1,x_2,\dots \} \, .
\end{equation}
The second marginal of $ P $ is called the distribution of the selector.
\end{definition}

Less formally, a selector is a randomized choice of a single
element. The conditional distribution $ P_x $ of $ z $ given $ x $ is a
probability measure concentrated on $ \{ x_1,x_2,\dots \} $. This
measure may happen to be a single atom, which leads to a non-randomized
selector
\begin{equation}\label{3n3}
z = x_{N(x_1,x_2,\dots)} \, ,
\end{equation}
where $ N $ is a Borel map $ (0,1)^\infty \to \{ 1,2,\dots \} $. (See
also \ref{3n9}.)

In order to prove existence of selectors with prescribed distributions
we use a deep duality theory for measures with given marginals, due to
Kellerer. It holds for a wide class of spaces $ \X_1, \X_2 $, but we
need only two. Below, in \ref{2n4}, \ref{2n5} and \ref{2n6} we assume
that
\begin{align*}
\X_1 \text{ is either } (0,1) \text{ or } (0,1)^\infty \, , \\
\X_2 \text{ is either } (0,1) \text{ or } (0,1)^\infty \, .
\end{align*}
Here is the result used here and once again in Sect.~\ref{sect8}.

\begin{theorem}\label{2n4} (Kellerer)
Let $ \mu_1, \mu_2 $ be probability measures on $ \X_1, \X_2 $
respectively, and $ B \subset \X_1 \times \X_2 $ a Borel set. Then
\[
S_{\mu_1,\mu_2}(B) = I_{\mu_1,\mu_2}(B) \, ,
\]
where $ S_{\mu_1,\mu_2}(B) $ is the supremum of $ \mu(B) $ over all probability
measures $ \mu $ on $ \X_1 \times \X_2 $ with marginals $ \mu_1, \mu_2
$, and $ I_{\mu_1,\mu_2}(B) $ is the infimum of $ \mu_1(B_1) +
\mu_2(B_2) $ over all Borel sets $ B_1 \subset \X_1 $, $ B_2 \subset
\X_2 $ such that $ B \subset (B_1\times\X_2) \cup (\X_1\times B_2) $.
\end{theorem}

See \cite[Corollary 2.18 and Proposition 3.3]{Ke84}.

Note that $ B $ need not be closed.

In fact, the infimum $ I_{\mu_1,\mu_2}(B) $ is always reached
\cite[Prop.~3.5]{Ke84}, but the supremum $ S_{\mu_1,\mu_2}(B) $ is not
always reached \cite[Example 2.20]{Ke84}.

\begin{remark}\label{2n5}
(a) If $ \mu_1,\mu_2 $ are positive (not just probability) measures such
that $ \mu_1(\X_1) = \mu_2(\X_2) $ then still $ S_{\mu_1,\mu_2}(B) =
I_{\mu_1,\mu_2}(B) $.

(b) $ S_{\mu_1,\mu_2}(B) $ is equal to the supremum of $ \nu(B) $ over
all positive (not just probability) measures $ \nu $ on $ B $ such that
$ \nu_1 \le \mu_1 $ and $ \nu_2 \le \mu_2 $, where $ \nu_1,\nu_2 $ are
the marginals of $ \nu $. This new supremum in $ \nu $ is reached if and
only if the original supremum in $ \mu $ is reached.
\end{remark}

\begin{lemma}\label{2n6}
Let $ \mu_1, \mu_2 $ be probability measures on $ \X_1, \X_2 $
respectively and $ B \subset \X_1 \times \X_2 $ a Borel set such that
\[
S_{\mu_1-\nu_1,\mu_2-\nu_2}(B) = S_{\mu_1,\mu_2}(B) - \nu(B)
\]
for every positive measure $ \nu $ on $ B $ such that $ \nu_1 \le
\mu_1 $ and $ \nu_2 \le \mu_2 $, where $ \nu_1,\nu_2 $ are the marginals
of $ \nu $. Then the supremum $ S_{\mu_1,\mu_2}(B) $ is reached.
\end{lemma}

(See also \ref{2.9}.)

\begin{proof}
First, taking a positive measure $ \nu $ on $ B $ such that $ \nu_1 \le
\mu_1 $, $ \nu_2 \le \mu_2 $ and $ \nu(B) \ge \frac12 S_{\mu_1,\mu_2}(B)
$ we get
\[
S_{\mu_1-\nu_1,\mu_2-\nu_2}(B) = S_{\mu_1,\mu_2}(B) - \nu(B) \le \frac12
S_{\mu_1,\mu_2}(B) \, .
\]
Second, taking a positive measure $ \nu' $ on $ B $ such that $ \nu'_1
\le \mu_1 - \nu_1 $, $ \nu'_2 \le \mu_2 - \nu_2 $ and $ \nu'(B) \ge
\frac12 S_{\mu_1-\nu_1,\mu_2-\nu_2}(B) $ we get $ (\nu+\nu')_1 \le \mu_1
$ and $ (\nu+\nu')_2 \le \mu_2 $, thus,
\[
S_{\mu_1-\nu_1-\nu'_1,\mu_2-\nu_2-\nu'_2}(B) = S_{\mu_1,\mu_2}(B) -
\nu(B) - \nu'(B) \le \frac14 S_{\mu_1,\mu_2}(B) \, .
\]
Continuing this way we get a convergent series of positive measures, $
\nu + \nu' + \nu'' + \dots $; its sum is a measure that reaches the
supremum indicated in \ref{2n5}(b).
\end{proof}

\begin{lemma}\label{3n7}
Let a random countable set $ \{ X_1, X_2, \dots \} $ satisfy
\begin{equation}\label{3n8}
\begin{gathered}
\text{for every Borel set $ B \subset (0,1) $ of positive measure,} \\
B \cap \{ X_1,X_2,\dots \} \ne \emptyset \text{ a.s.}
\end{gathered}
\end{equation}
Then the random set has a selector distributed uniformly on $ (0,1) $.
\end{lemma}

\begin{proof}
We apply Theorem \ref{2n4} to $ \X_1 = (0,1)^\infty $, $ \X_2 = (0,1) $,
$ \mu_1 $ --- the joint distribution of $ X_1,X_2,\dots $, $ \mu_2 $ ---
the uniform distribution on $ (0,1) $, and $ B $ --- the set of all
pairs $ (x,z) $ satisfying \eqref{3n2}. By \eqref{3n8}, $ B $ intersects
$ B_1 \times B_2 $ for all Borel sets $ B_1 \subset \X_1 $, $ B_2
\subset \X_2 $ of positive measure. Therefore $ I_{\mu_1,\mu_2}(B) = 1
$. By the same argument, all absolutely continuous measures $
\nu_1,\nu_2 $ on $ \X_1,\X_2 $ respectively, such that $ \nu_1 (\X_1) =
\nu_2 (\X_2) $, satisfy $ I_{\nu_1,\nu_2} (B) = \nu_1 (\X_1) $. By
\ref{2n5}(a), $ S_{\nu_1,\nu_2} (B) = I_{\nu_1,\nu_2} (B) $. Thus, the
condition of  Lemma \ref{2n6} is satisfied (by $ \mu_1,\mu_2 $ and $ B
$). By \ref{2n6}, some measure $ P $ reaches $ S_{\mu_1,\mu_2}(B) =
I_{\mu_1,\mu_2}(B) = 1 $ and therefore $ P $ is the needed selector.
\end{proof}

\begin{counterexample}\label{3n9}
In Lemma \ref{3n7} one cannot replace `a selector' with `a non-randomized
selector \eqref{3n3}'. Randomization is essential!

Let $ \Om = (0,1) $ (with Lebesgue measure) and $ \{ X_1(\om), X_2(\om),
\dots \} = ( \frac12 \om + \Q ) \cap (0,1) $ where $ \Q \subset \R $ is
the set of rational numbers (and $ \frac12 \om + \Q $ is its shift by $
\frac12 \om $). Then \eqref{3n8} is satisfied (since $ B + \Q $ is of full
measure), but every selector $ Z : \Om \to (0,1) $ of the form
$ Z(\om) = X_{N(\om)} (\om) $ has a nonuniform distribution. \emph{Proof:}
let $ A_q = \{ \om\in(0,1) : Z(\om) - \frac12 \om = q \} $ for $ q \in
\Q $, then
\[
\Pr{ Z \in B } = 2 \sum_{q\in\Q} \int_B \One_{A_q} (2(x-q)) \, \D x \, ,
\]
which shows that the distribution of $ Z $ has a density taking on the
values $ 0,2,4,\dots $ only.
\end{counterexample}

\section[]{\raggedright Independence}
\label{sect4}
According to \eqref{2n0}, our `random countable set' $ \{X_1,X_2,\dots\}
$ can be finite, but cannot be empty. This is why in general we cannot
treat the intersection, say, $ \{X_1,X_2,\dots\} \cap (0,\frac12) $ as a
random countable set.

By a \emph{random dense countable subset of $ (0,1) $} we mean a random
countable subset $ \{X_1,X_2,\dots\} $ of $ (0,1) $, dense in $ (0,1) $
a.s. Equivalently, $ \Pr{ \exists k \;\; a<X_k<b } = 1 $ whenever $ 0
\le a < b \le 1 $. A random dense countable subset of another interval
is defined similarly. It is easy to see that $ \{X_1,X_2,\dots\} \cap
(a,b) $ is a random dense countable subset of $ (a,b) $ whenever $
\{X_1,X_2,\dots\} $ is a random dense countable subset of $ (0,1) $ and
$ (a,b) \subset (0,1) $. We call $ \{X_1,X_2,\dots\} \cap (a,b) $ a
\emph{fragment} of $ \{X_1,X_2,\dots\} $.

It may happen that two (or more) fragments can be described by
independent sequences of random variables; such fragments will be called
independent. The definition is formulated below in terms of random
variables, but could be reformulated in terms of measures on $
(0,1)^\infty $.

\begin{definition}\label{4n1}
Let $ \{X_1,X_2,\dots\} $ be a random dense countable subset of $ (0,1)
$. We say that two fragments $ \{X_1,X_2,\dots\} \cap (0,\frac12) $ and
$ \{X_1,X_2,\dots\} \cap [\frac12,1) $ of $ \{X_1,X_2,\dots\} $ are
\emph{independent,} if there exist random variables $ Y_1,Y_2,\dots $
(on some probability space) such that

(a) $ \{ Y_1,Y_2,\dots \} $ is distributed like $ \{X_1,X_2,\dots\} $;

(b) $ Y_{2k-1} < \frac12 $ and $ Y_{2k} \ge \frac12 $ a.s.\ for $
k=1,2,\dots $;

(c) the random sequence $ (Y_2,Y_4,Y_6,\dots) $ is independent of the
random sequence $ (Y_1,Y_3,Y_5,\dots) $.

Similarly we define independence of $ n $ fragments $
\{X_1,X_2,\dots\} \cap [a_{k-1},a_k) $, $ k=1,\dots,n $, for any $
n=2,3,\dots $ and any $ a_0,\dots,a_n $ such that $
0=a_0<a_1<\dots<a_n=1 $.
\end{definition}

\begin{definition}\label{4n2}
A random dense countable subset $ \{X_1,X_2,\dots\} $ of $ (0,1) $
satisfies \emph{the independence condition,} if for every $ n=2,3,\dots $
and every $ a_0,\dots,a_n $ such that $ 0=a_0<a_1<\dots<a_n=1 $ the $ n
$ fragments $ \{X_1,X_2,\dots\} \cap [a_{k-1},a_k) $ ($ k=1,\dots,n $)
are independent.
\end{definition}

Such random dense countable sets are described below, assuming that each
$ X_k $ has a density, in other words,
\begin{equation}\label{4n3}
\begin{gathered}
\text{for every Borel set $ B \subset (0,1) $ of measure $ 0 $,} \\
B \cap \{ X_1,X_2,\dots \} = \emptyset \quad \text{a.s.}
\end{gathered}
\end{equation}
(In contrast to Main lemma, existence of \emph{joint} densities is not
assumed. See also \ref{5n10}.)

\begin{proposition}\label{4n4}
For every random dense countable set satisfying \eqref{4n3} and the
independence condition there exists a measurable function $ r : (0,1)
\to [0,\infty] $ such that for every Borel set $ B \subset (0,1) $

(a) if $ \int_B r(x) \, \D x = \infty $ then the set $ B \cap \{
X_1,X_2,\dots \} $ is infinite a.s.;

(b) if $ \int_B r(x) \, \D x < \infty $ then the set $ B \cap \{
X_1,X_2,\dots \} $ is finite a.s., and the number of its elements has
the Poisson distribution with the mean $ \int_B r(x) \, \D x $.
\end{proposition}

Note that $ r(\cdot) $ may be infinite. See also \ref{9n4}.

The proof is given after some remarks and lemmas.

\begin{remark}\label{4n5}
The function $ r $ is determined uniquely (up to equality almost
everywhere) by the random dense countable set and moreover, by its
distribution. That is, if $ \{ Y_1,Y_2,\dots \} $ is distributed like $
\{X_1,X_2,\dots\} $ and \ref{4n4}(a,b) hold for $ \{ Y_1,Y_2,\dots \} $
and another function $ r_1 $, then $ r_1(x) = r(x) $ for almost all $ x
\in (0,1) $.
\end{remark}

\begin{remark}
The function $ r $ is just the sum
\[
r(\cdot) = f_1(\cdot) + f_2(\cdot) + \dots
\]
of the densities $ f_k $ of $ X_k $.
\end{remark}

\begin{remark}
For every measurable function $ r : (0,1) \to [0,\infty] $ such that $
\int_a^b r(x) \, \D x = \infty $ whenever $ 0 \le a < b \le 1 $, there
exists a random dense countable set $ \{ X_1,X_2,\dots \} $ satisfying
\eqref{4n3}, the independence condition and \ref{4n4}(a,b) (for the
given $ r $). See also \ref{9n5} and \cite[Sect.~2.5]{Ki}.

If $ r(x)=\infty $ almost everywhere, we take an unordered infinite
sample.

If $ r(x)<\infty $ almost everywhere, we take a Poisson point process
with the intensity measure $ r(x) \, \D x $.

Otherwise we combine an unordered infinite sample on $ \{ x :
r(x)=\infty \} $ and a Poisson point process on $ \{ x : r(x)<\infty \}
$.
\end{remark}

In order to prove \ref{4n4}, for a given $ B \subset (0,1) $ we denote
by $ \xi(x,y) $ the (random) number of elements (maybe, $ \infty $) in
the set $ B \cap (x,y) \cap \{X_1,X_2,\dots\} $ and introduce
\begin{align*}
\al(x,y) &= \Pr{ \xi(x,y) = 0 } \, ,
 \\
\be(x,y) &= \Ex \exp ( - \xi(x,y) )
\end{align*}
for $ 0 \le x < y \le 1 $. (Of course, $ \exp(-\infty) = 0 $.) Clearly,
\begin{equation}\label{4n7}
\( 1 - \E^{-1} \) \( 1 - \al(x,y) \) \le 1 - \be(x,y) \le 1 -
\al(x,y)
\end{equation}
for $ 0 \le x < y \le 1 $. By \eqref{4n3} and the independence
condition,
\begin{align*}
\xi(x,y) + \xi(y,z) &= \xi(x,z) \, , \\
\al(x,y) \al(y,z) &= \al(x,z) \, , \\
\be(x,y) \be(y,z) &= \be(x,z)
\end{align*}
whenever $ 0 \le x < y < z \le 1 $.

\begin{lemma}\label{4n8}
If $ \be(0,1) \ne 0 $ then $ \be(x-\eps,x+\eps) \to 1 $ as $ \eps\to0+ $
for every $ x \in (0,1) $.
\end{lemma}

\begin{proof}
Let $ x_1 < x_2 < \dots $, $ x_k \to x $. Random variables $
\xi(x_k,x_{k+1}) $ are independent. By Kolmogorov's \OI\ law, the event
$ \xi(x_k,x) \to 0 $ is of probability $ 0 $ or $ 1 $. It cannot be of
probability $ 0 $, since then $ \xi(x_1,x) = \infty $ a.s., which
implies $ \be(0,1) = 0 $. Thus, $ \xi(x_k,x) \to 0 $ a.s., therefore $
\be(x_k,x) \to 1 $ and $ \be(x-\eps,1) \to 1 $. Similarly, $
\be(x,x+\eps) \to 1 $.
\end{proof}

\begin{lemma}\label{4n9}
If $ \be(0,1) \ne 0 $ then $ \al(0,1) \ne 0 $.
\end{lemma}

\begin{proof}
By \ref{4n8}, for every $ x \in (0,1) $ there exists $ \eps > 0 $ such
that $ \be(x-\eps,x+\eps) > \E^{-1} $, therefore $ \al(x-\eps,x+\eps)
\ne 0 $ by \eqref{4n7}. (For $ x=0 $, $ x=1 $ we use one-sided
neighborhoods.) Choosing a finite covering and using multiplicativity of
$ \al $ we get $ \al(0,1) \ne 0 $.
\end{proof}

\begin{lemma}\label{4n10}
If $ \al(0,1) \ne 0 $ then $ \xi(0,1) $ has the Poisson distribution
with the mean $ -\ln \al(0,1) $.
\end{lemma}

\begin{proof}
By \ref{4n8} and \eqref{4n7}, $ \al(x-\eps,x+\eps) \to 1 $. We define a
nonatomic finite positive measure $ \mu $ on $ [0,1] $ by
\[
\mu([x,y]) = -\ln \al(x,y) \quad \text{for } 0 \le x < y \le 1
\]
and introduce a Poisson point process on $ [0,1] $ whose intensity
measure is $ \mu $. Denote by $ \eta(x,y) $ the (random) number of
Poisson points on $ [x,y] $, then $ \Ex \eta(x,y) = \mu([x,y]) $ and $
\Pr{ \eta(x,y)=0 } = \exp \( -\Ex \eta(x,y) \) = \al(x,y) = \Pr{
\xi(x,y)=0 } $. In other words, the two random variables $ \xi(x,y) \wedge
1 $ and $ \eta(x,y) \wedge 1 $ are identically distributed (of course, $ a
\wedge b = \min(a,b) $). By independence, for any $ n $ the joint
distribution of $ n $ random variables $ \xi\(\frac{k-1}n,\frac k n\)
\wedge 1 $ (for $ k=1,\dots,n $) is equal to the joint distribution of $ n
$ random variables $ \eta\(\frac{k-1}n,\frac k n\) \wedge 1 $. Taking
into account that
\[
\xi(0,1) = \lim_{n\to\infty} \sum_{k=1}^n \bigg(
\xi\Big(\frac{k-1}n,\frac k n\Big) \wedge 1 \bigg)
\]
and the same for $ \eta $, we conclude that $ \xi(0,1) $ and $ \eta(0,1)
$ are identically distributed.
\end{proof}

\begin{remark}
In addition (but we do not need it),

(a) the joint distribution of $ \xi(r,s) $ for all rational $ r,s $ such
that $ 0 \le r < s \le 1 $ (this is a countable family of random
variables) is equal to the joint distribution of all $ \eta(r,s) $,

(b) the random finite set $ B \cap \{X_1,X_2,\dots\} $ is distributed
like the Poisson point process,

(c) the measure $ \mu $ has the density $ (f_1+f_2+\dots ) \cdot \One_B
$, where $ f_k $ is the density of $ X_k $.
\end{remark}

\begin{proof}[Proof of Proposition \textup{\ref{4n4}}]
We take $ r = f_1 + f_2 + \dots $, note that $ \int_B r(x) \, \D x = \Ex
\xi(0,1) $ and prove (b) first.

(b) Let $ \int_B r(x) \, \D x < \infty $, then $ \xi(0,1) < \infty $
a.s., therefore $ \be(0,1) \ne 0 $. By \ref{4n9}, $ \al(0,1) \ne 0 $. By
\ref{4n10}, $ \xi(0,1) $ has a Poisson distribution.

(a) Let $ \int_B r(x) \, \D x = \infty $, then $ \Ex \xi(0,1) = \infty
$, thus $ \xi(0,1) $ cannot have a Poisson distribution. By the argument
used in the proof of (b), $ \be(0,1) = 0 $. Therefore $ \xi(0,1) =
\infty $ a.s.
\end{proof}

\section[]{\raggedright Selectors and independence}
\label{sect5}We consider a random dense countable subset $ \{X_1,X_2,\dots\} $ of $
(0,1) $, satisfying the independence condition and \eqref{3n8}. If
\eqref{4n3} is also satisfied then \eqref{3n8} means that the
corresponding function $ r $ (see \ref{4n4}) is infinite almost
everywhere.

A uniformly distributed selector exists by \ref{3n7}. Moreover, there
exists a pair of independent uniformly distributed selectors. It follows
via Th.~\ref{2n4} from the fact that $ \{X_1,X_2,\dots\} \times
\{X_1,X_2,\dots\} $ intersects a.s.\ any given set $ B \subset (0,1)
\times (0,1) $ of positive measure. \emph{Hint:} we may assume that $ B
\subset (0,\theta) \times (\theta,1) $ for some $ \theta \in (0,1) $;
consider independent fragments $ \{ Y_1,Y_3,\dots \} = (0,\theta) \cap
\{X_1,X_2,\dots\} $, $ \{ Y_2,Y_4,\dots \} = [\theta,1) \cap
\{X_1,X_2,\dots\} $; almost surely, the first fragment intersects the
first projection of $ B $, and the second fragment intersects the
corresponding section of $ B $.

However, we need a stronger statement: for every selector $ Z_1 $ there
exists a selector $ Z_2 $ distributed uniformly and independent of $ Z_1
$; here is the exact formulation. (The proof is given after Lemma
\ref{5n8}.)

\begin{proposition}\label{5n1}
Let $ \{X_1,X_2,\dots\} $ be a random dense countable subset of $ (0,1)
$ satisfying the independence condition and \eqref{3n8}. Let a
probability measure $ P_1 $ on $ (0,1)^\infty \times (0,1) $ be a
selector of $ \{X_1,X_2,\dots\} $ (as defined by \ref{3n1}). Then there
exists a probability measure $ P_2 $ on $ (0,1)^\infty \times (0,1)^2 $
such that, denoting points of $ (0,1)^\infty \times (0,1)^2 $ by $
(x,(z_1,z_2)) $, we have (w.r.t.\ $ P_2 $)

(a) the joint distribution of $ x $ and $ z_1 $ is equal to $ P_1 $;

(b) the distribution of $ z_2 $ is uniform on $ (0,1) $;

(c) $ z_1, z_2 $ are independent;

(d) $ z_2 \in \{ x_1,x_2,\dots \} $ a.s. (where $ (x_1,x_2,\dots) = x $).
\end{proposition}

Conditioning on $ z_1 $ decomposes the two-selector problem into a
continuum of single-selector problems. In terms of conditional
distributions $ P_1 ( \D x | z_1 ) $, $ P_2 ( \D x \D z_2 | z_1 ) $ we
need the following:

(e) $ z_2 \in \{ x_1,x_2,\dots \} $ for $ P_2(\cdot|z_1) $-almost all
$ ((x_1,x_2,\dots),z_2) $;

(f) the distribution of $ x $ according to $ P_2 ( \D x \D z_2 | z_1 ) $
is $ P_1 ( \D x | z_1 ) $;

(g) the distribution of $ z_2 $ according to $ P_2 ( \D x \D z_2 | z_1 )
$ is uniform on $ (0,1) $.

\noindent That is, we need a uniformly distributed selector of a random
set distributed $ P_1 ( \cdot | z_1 ) $. To this end we will transfer
\eqref{3n8} from the unconditional joint distribution of $ X_1,X_2,\dots
$ to their conditional joint distribution $ P_1 ( \cdot | z_1 ) $.

\begin{proposition}\label{5n2}
Let $ X_1,X_2,\dots $ and $ Y_1,Y_2,\dots $ be as in Def.~\ref{4n1}, and
$ \Pr{ X_1 < \frac12 } > 0 $. Then for almost all $ x_1 \in (0,\frac12)
$ (w.r.t.\ the distribution of $ X_1 $), the conditional joint
distribution of $ Y_2,Y_4,\dots $ given $ X_1=x_1 $ is absolutely
continuous w.r.t.\ the (unconditional) joint distribution of $
Y_2,Y_4,\dots $ 
\end{proposition}

The proof is given before Lemma \ref{5n7}.

\begin{counterexample}\label{5n3}
Condition \ref{4n1}(c) is essential for \ref{5n2} (in spite of the fact
that $ Y_1,Y_3,\dots $ are irrelevant).

We take independent random variables $ Z_1,Z_2,\dots $ such that each $
Z_{2k-1} $ is uniform on $ (0,\frac12) $ and each $ Z_{2k} $ is uniform on
$ (\frac12,1) $. We define $ Y_1,Y_2,\dots $ as follows. First, $
Y_{2k-1} = Z_{2k-1} $. Second, if the $ k $-th binary digit of $ 2Y_1 $
is equal to $ 1 $ then $ Y_{4k-2} = \min(Z_{4k-2},Z_{4k}) $ and $
Y_{4k} = \max(Z_{4k-2},Z_{4k}) $; otherwise (if the digit is $ 0 $), $
Y_{4k-2} = \max(Z_{4k-2},Z_{4k}) $ and $ Y_{4k} = \min(Z_{4k-2},Z_{4k})
$.

\begin{sloppypar}
We get a random dense countable set $ \{Y_1,Y_2,\dots\} $ whose
fragments $ \{Y_1,Y_2,\dots\} \cap (0,\frac12) = \{Y_1,Y_3,\dots\} $ and
$ \{Y_1,Y_2,\dots\} \cap (\frac12,1) = \{Y_2,Y_4,\dots\} $ are
independent. However, $ Y_1 $ is a function of $ Y_2,Y_4,\dots $ (and of
course, the conditional joint distribution of $ Y_2,Y_4,\dots $ given $
Y_1 $ is singular to their unconditional joint distribution).
\end{sloppypar}
\end{counterexample}

In order to prove \ref{5n2} we may partition the event $ X_1 < \frac12 $
into events $ X_1 = Y_{2k-1} $. Within such event the condition $ X_1 =
x_1 $ becomes just $ Y_{2k-1} = x_1 $. However, it does not make the
matter trivial, since the event $ X_1 = Y_{2k-1} $ need not belong to
the \sif\ generated by $ Y_1,Y_3,\dots $ (nor to the \sif\ generated by
$ X_1 $).

\smallskip
\textsc{digression: nonsingular pairs}
\smallskip

Sometimes dependence between two random variables reduces to a joint density
(w.r.t.\ their marginal distributions). Here are two formulation in general
terms.

\begin{lemma}\label{5n4}
Let $ (\Om,\F,P) $ be a probability space and $ C \subset \Om $ a measurable
set. The following two conditions on a pair of sub-\sif s $ \F_1, \F_2 \subset
\F $ are equivalent:

(a) there exists a measurable function $ f : \Om \times \Om \to [0,\infty) $
such that
\[
P ( A \cap B \cap C ) = \int_{A\times B} f(\om_1,\om_2) \, P(\D\om_1)
P(\D\om_2) \quad \text{for all } A \in \F_1, B \in \F_2 \, ;
\]

(b) there exists a measurable function $ g : C \times C \to [0,\infty) $
such that
\[
P ( A \cap B \cap C ) = \int_{(A\cap C)\times(B\cap C)} g(\om_1,\om_2) \,
P(\D\om_1) P(\D\om_2) \quad \text{for all } A \in \F_1, B \in \F_2 \, .
\]
\end{lemma}

(Note that $ f,g $ may vanish somewhere, and $ C $ need not belong to $
\F_1 $ or $ \F_2 $.)

\begin{proof}
(b) \imp (a): just take $ f(\om_1,\om_2) = g(\om_1,\om_2) $ for $ \om_1,\om_2
  \in C $ and $ 0 $ otherwise.

(a) \imp (b):
we consider conditional probabilities $ h_1 = \cP{C}{\F_1} $, $ h_2 =
\cP{C}{\F_2} $, note that $ h_1(\om) > 0 $, $ h_2(\om) > 0 $ for almost all $
\om \in C $ and define
\[
g(\om_1,\om_2) = \frac{ f(\om_1,\om_2) }{ h_1(\om_1) h_2(\om_2) } \quad
\text{for } \om_1, \om_2 \in C \, .
\]
Then
\begin{multline*}
\int_{(A\cap C)\times(B\cap C)} g(\om_1,\om_2) \, P(\D\om_1) P(\D\om_2) = \\
= \int_{A\times B} \frac{ f(\om_1,\om_2) \One_C(\om_1) \One_C(\om_2) }{
 h_1(\om_1) h_2(\om_2) } \, P(\D\om_1) P(\D\om_2) \, .
\end{multline*}
(The integrand is treated as $ 0 $ outside $ C \times C $.)
Assuming that $ f $ is \measurable{(\F_1\otimes\F_2)} (otherwise $ f $ may be
replaced with its conditional expectation) we see that the conditional
expectation of the integrand, given $ \F_1\otimes\F_2 $, is equal to $
f(\om_1,\om_2) $. Thus, the integral is
\[
\dots = \int_{A\times B} f(\om_1,\om_2) \, P(\D\om_1) P(\D\om_2) = P ( A \cap
B \cap C ) \, .
\]
\end{proof}

\begin{definition}\label{5n5}
Let $ (\Om,\F,P) $ be a probability space and $ C \subset \Om $ a measurable
set. Two sub-\sif s $ \F_1, \F_2 \subset \F $ are a \emph{nonsingular pair}
within $ C $, if they satisfy the equivalent conditions of Lemma \ref{5n4}.
\end{definition}

\begin{lemma}\label{5n6}
(a)
Let $ C_1 \subset C_2 $. If $ \F_1, \F_2 $ are a nonsingular pair within $ C_2
$ then they are a nonsingular pair within $ C_1 $.

(b)
Let $ C_1,C_2,\dots $ be pairwise disjoint and $ C = C_1 \cup C_2 \cup \dots $
If $ \F_1, \F_2 $ are a nonsingular pair within $ C_k $ for each $ k $ then
they are a nonsingular pair within $ C $.

(c)
Let $ \Ec_1 \subset \F $ be another sub-\sif\ such that $ \Ec_1 \subset \F_1 $
within $ C $ in the sense that
\[
\forall E \in \Ec_1 \;\; \exists A \in \F_1 \;\; ( A \cap C = E \cap C ) \, .
\]
If $ \F_1, \F_2 $ are a nonsingular pair within $ C $ then $ \Ec_1, \F_2 $
are a nonsingular pair within $ C $.
\end{lemma}

\begin{proof}
(a)
We define two measures $ \mu_1, \mu_2 $ on $ (\Om,\F_1) \times (\Om,\F_2) $ by
$ \mu_k (Z) = P ( C_k \cap \{ \om : (\om,\om) \in Z \} ) $ for $ k=1,2
$. Clearly, $ \mu_k (A\times B) = P ( A \cap B \cap C_k ) $. Condition
\ref{5n4}(a) for $ C_k $ means absolute continuity of $ \mu_k $ (w.r.t.\
$ P |_{\F_1} \times P |_{\F_2} $). However, $ \mu_1 \le \mu_2 $.

(b)
Using the first definition, \ref{5n4}(a), we just take $ f = f_1 + f_2 + \dots
$

(c)
Immediate, provided that the second definition us used, \ref{5n4}(b).
\end{proof}

\textsc{end of digression}
\smallskip

\begin{proof}[Proof of Prop.~\textup{\ref{5n2}}]
It is sufficient to prove that the two \sif s $ \si(X_1) $ (generated by
$ X_1 $) and $ \si(Y_2,Y_4,\dots) $ are a nonsingular pair within the
event $ C = \{ X_1 < \frac12 \} $. Without loss of generality we assume
that $ Y_1, Y_3, \dots $ are pairwise different a.s. (otherwise we skip
redundant elements via a random renumbering). We partition $ C $ into
events $ C_k = \{ X_1 = Y_{2k-1} \} $. Lemma \ref{5n6}(b) reduces $ C $
to $ C_k $. By \ref{5n6}(c) we replace $ \si(X_1) $ with $ \si(Y_{2k-1})
$. By \ref{5n6}(a) we replace $ C_k $ with the whole $ \Om $. Finally,
the \sif s $ \si(Y_{2k-1}) $, $ \si(Y_2,Y_4,\dots) $ are a nonsingular
pair within $ \Om $, since they are independent.
\end{proof}

\begin{lemma}\label{5n7}
Let $ \{X_1,X_2,\dots\} $ be a random dense countable subset of $ (0,1)
$ satisfying the independence condition and \eqref{3n8}. Then for every
Borel set $ B \subset (0,1) $ of positive measure,
\[
\cP{ \{ X_2,X_3,\dots \} \cap B \ne \emptyset }{ X_1 } = 1 \quad
\text{a.s.}
\]
\end{lemma}

\begin{proof}
First, we assume in addition that $ \mes \( B \cap (\frac12,1) \) > 0 $
(`$ \mes $' stands for Lebesgue measure) and $ \Pr{ X_1 < \frac12 } = 1
$. Introducing $ Y_k $ according to Def.~\ref{4n1} we note that $ \{
Y_2,Y_4,\dots \} \cap B \supset \{ X_1,X_2,\dots \} \cap B \cap
(\frac12,1) \ne \emptyset $ a.s.\ by \eqref{3n8}. It follows via
Prop.~\ref{5n2} that $ \cP{ \{ Y_2,Y_4,\dots \} \cap B \ne \emptyset }{
X_1 } = 1 $ a.s. Taking into account that $ \{ X_2,X_3,\dots \} \supset
\{ Y_2,Y_4,\dots \} $ we get $ \cP{ \{ X_2,X_3,\dots \} \cap B \ne
\emptyset }{ X_1 } = 1 $ a.s.

Similarly we consider the case $ \mes \( B \cap (0,\frac12) \) > 0 $
and $ \Pr{ X_1 \ge \frac12 } = 1 $.

Assuming both $ \mes \( B \cap (0,\frac12) \) > 0 $ and $ \mes \( B \cap
(\frac12,1) \) > 0 $ we get the same conclusion for arbitrary
distribution of $ X_1 $.

The same arguments work for any threshold $ \theta \in (0,1) $ instead
of $ \frac12 $. It remains to note that for every $ B $ there exists $
\theta $ such that both $ \mes \( B \cap (0,\theta) \) > 0 $ and $ \mes
\( B \cap (\theta,1) \) > 0 $.
\end{proof}

The claim of Lemma \ref{5n7} is of the form $ \forall B \; \(
\cP{\dots}{X_1} = 1 \text{ a.s.} \) $, but the following lemma gives
more: $ \cP{ \forall B \; (\dots) }{X_1} = 1 $ a.s.

\begin{lemma}\label{5n8}
Let $ \{X_1,X_2,\dots\} $ be a random dense countable subset of $ (0,1)
$ satisfying the independence condition and \eqref{3n8}. Denote by $ \nu
$ the distribution of $ X_1 $ and by $ \mu_{x_1} $ the conditional joint
distribution of $ X_2,X_3,\dots $ given $ X_1 = x_1 $. (Of course, $
\mu_{x_1} $ is well-defined for \almost{\nu} all $ x_1 $.) Then
\almost{\nu} all $ x_1 \in (0,1) $ are such that for every Borel set $ B
\subset (0,1) $ of positive measure,
\[
\mu_{x_1} \( \big\{ (x_2,x_3,\dots) : \{x_2,x_3,\dots\} \cap B \ne
\emptyset \big\} \) = 1 \, .
\]
\end{lemma}

\begin{proof}
The proof of \ref{5n7} needs only tiny modification, but the last
paragraph (about $ \theta $) needs some attention. The exceptional set
of $ x_1 $ may depend on $ \theta $, which is not an obstacle since we
may use only rational $ \theta $. Details are left to the reader.
\end{proof}

\begin{proof}[Proof of Prop.~\textup{\ref{5n1}}]
We apply Lemma \ref{5n8} to the sequence $ (Z_1,X_1,X_2,\dots) $ rather
than $ (X_1,X_2,\dots) $; here $ Z_1 $ is the given selector. More
formally, we consider the image of the given measure $ P_1 $ under the
map $ (0,1)^\infty \times (0,1) \to (0,1)^\infty $ defined by $ \(
(x_1,x_2,\dots),z_1 \) \mapsto ( z_1, x_1,x_2,\dots ) $.

Lemma \ref{5n8} introduces $ \nu $ (the distribution of $ Z_1 $) and $
\mu_{z_1} $ (the conditional distribution of $ (X_1,X_2,\dots) $ given $
Z_1=z_1 $), and states that \eqref{3n8} is satisfied by $ \mu_{z_1} $
for \almost{\nu} all $ z_1 $. Applying \ref{3n7} to $ \mu_{z_1} $ we get
a probability measure $ \ti\mu_{z_1} $ on $ (0,1)^\infty \times (0,1) $
such that the first marginal of $ \ti\mu_{z_1} $ is equal to $ \mu_{z_1}
$, the second marginal of $ \ti\mu_{z_1} $ is the uniform distribution
on $ (0,1) $, and \almost{\ti\mu_{z_1}} all pairs $ \(
(x_1,x_2,\dots),z_2 \) $ satisfy $ z_2 \in \{ x_1,x_2,\dots \} $.

In order to combine measures $ \ti\mu_{z_1} $ into a measure $ P_2
$ we need measurability of the map $ z_1 \mapsto \ti\mu_{z_1} $.

The set of all probability measures on $ (0,1)^\infty \times (0,1) $ is
a standard Borel space (see \cite{Ke95}, Th.~(17.24) and the paragraph
after it), and the map $ \mu \mapsto \mu(B) $ is Borel for every Borel
set $ B \subset (0,1)^\infty \times (0,1) $. (In fact, these maps
generate the Borel \sif\ on the space of measures.) It follows easily
that the subset $ M $ of the space of measures, introduced below, is
Borel. Namely, $ M $ is the set of all $ \mu $ such that the second
marginal of $ \mu $ is the uniform distribution on $ (0,1) $ and $ \mu $
is concentrated on the set of $ \( (x_1,x_2,\dots), z_2 \) $ such that $
z_2 \in \{ x_1,x_2,\dots\} $. Also, the first marginal of $ \mu $ is a
Borel function of $ \mu $ (which means a Borel map from the space of
measures on $ (0,1)^\infty \times (0,1) $ into the similar space of
measures on $ (0,1)^\infty $).

The conditional measure $ \mu_{z_1} $ is a \measurable{\nu} function of
$ z_1 $ defined \almost{\nu} everywhere; it may be chosen to be a Borel
map from $ (0,1) $ to the space of measures on $ (0,1)^\infty $. In
addition we can ensure that each $ \mu_{z_1} $ is the first marginal of
some $ \ti\mu_{z_1} \in M $. It follows that these $ \ti\mu_{z_1} \in M
$ can be chosen as a \measurable{\nu} (maybe not Borel, see
\cite[5.1.7]{Sri}) function of $ z_1 $, by the (Jankov and) von Neumann
uniformization theorem, see \cite[Sect.~18A]{Ke95} or
\cite[Sect.~5.5]{Sri}.

Now we combine these $ \ti\mu_{z_1} $ into a probability measure $ P_2 $
on $ (0,1)^\infty \times (0,1)^2 $ such that, denoting a point of $
(0,1)^\infty \times (0,1)^2 $ by $ \( (x_1,x_2,\dots), (z_1,z_2) \) $ we
have: $ z_1 $ is distributed $ \nu $, and $ P_2 ( \D x \D z_2 | z_1 ) =
\ti\mu_{z_1} (\D x \D z_2) $.

It remains to note that $ P_2 $ satisfies (e), (f), (g) formulated after
Prop.~\ref{5n1}. The first marginal of $ \ti\mu_{z_1} = P_2(\cdot|z_1) $
is equal to $ \mu_{z_1} = P_1(\cdot|z_1) $, which verifies (f). The
second marginal of $ \ti\mu_{z_1} = P_2(\cdot|z_1) $ is the uniform
distribution on $ (0,1) $, which verifies (g). And $ z_2 \in \{
x_1,x_2,\dots \} $ almost sure w.r.t.\ $ \ti\mu_{z_1} = P_2(\cdot|z_1)
$, which verifies (e).
\end{proof}

Prop.~\ref{5n1} is a special case ($ n=1 $) of Prop.~\ref{5n9} below;
the latter shows that for every $ n $ selectors $ Z_1,\dots,Z_n $ there
exists a selector $ Z_{n+1} $ distributed uniformly and independent of $
Z_1,\dots,Z_n $.

\begin{proposition}\label{5n9}
Let $ \{X_1,X_2,\dots\} $ be a random dense countable subset of $ (0,1)
$ satisfying the independence condition and \eqref{3n8}. Let $ n \in
\{1,2,\dots\} $ be given, and $ P_n $ be a probability measure on $
(0,1)^\infty \times (0,1)^n $ such that

(i) the first marginal of $ P_n $ is equal to the joint distribution of
$ X_1,X_2,\dots $;

(ii) \almost{P_n} all pairs $ (x,z) $, $ x=(x_1,x_2,\dots) $, $
z=(z_1,\dots,z_n) $ are such that $ \{z_1,\dots,z_n\} \subset
\{x_1,x_2,\dots\} $.

Then there exists a probability measure $ P_{n+1} $ on $ (0,1)^\infty
\times (0,1)^{n+1} $ such that, denoting points of $ (0,1)^\infty \times
(0,1)^{n+1} $ by $ (x,(z_1,\dots,z_{n+1})) $, we have (w.r.t.\ $ P_{n+1}
$)

(a) the joint distribution of $ x $ and $ (z_1,\dots,z_n) $ is equal to
$ P_n $;

(b) the distribution of $ z_{n+1} $ is uniform on $ (0,1) $;

(c) $ z_{n+1} $ is independent of $ (z_1,\dots,z_n) $;

(d) $ z_{n+1} \in \{ x_1,x_2,\dots \} $ a.s. (where $ (x_1,x_2,\dots) = x $).
\end{proposition}

The proof, quite similar to the proof of Prop.~\ref{5n1}, is left to the
reader, but some hints follow. Two independent fragments of a random set
are used in \ref{5n8}, according to the partition of $ (0,1) $ into $
(0,\theta) $ and $ [\theta,1) $, where $ \theta \in (0,1) $ is
rational. One part contains $ z_1 $, the other part contains a portion
of the given set $ B $ of positive measure. Now, dealing with $
z_1,\dots,z_n $ we still partition $ (0,1) $ in two parts, but they are
not just intervals. Rather, each part consists of finitely many
intervals with rational endpoints. Still, the independence condition
gives us two independent fragments.

Here is another implication of the independence condition. In some sense
the proof below is similar to the proof of \ref{5n1}, \ref{5n9}. There,
\eqref{3n8} was transferred to conditional distributions via
\ref{5n2}. Here we do it with \eqref{4n3}.

\begin{lemma}\label{5n10}
Let $ \{X_1,X_2,\dots\} $ be a random dense countable subset of $ (0,1)
$ satisfying the independence condition and \eqref{4n3}. If $ \Pr{ X_k =
X_l } = 0 $ whenever $ k \ne l $ then for every $ n $ the joint
distribution of $ X_1,\dots,X_n $ is absolutely continuous.
\end{lemma}

\begin{proof}
Once again, I restrict myself to the case $ n=2 $, leaving the general
case to the reader.

The marginal (one-dimensional) distribution of any $ X_n $ is absolutely
continuous due to \eqref{4n3}. It is sufficient to prove that the
conditional distribution of $ X_2 $ given $ X_1 $ is absolutely
continuous, that is, $ \cP{ X_2 \in B }{ X_1 } = 0 $ a.s.\ for all
negligible $ B \subset (0,1) $ simultaneously. By Prop.~\ref{5n2} it
holds for $ X_1 < \frac12 $ and $ B \subset (\frac12,1) $. Similarly, it
holds for $ X_1 < \theta $ and $ B \subset (\theta,1) $, or $ X_1 >
\theta $ and $ B \subset (0,\theta) $, for all rational $ \theta $
simultaneously. Therefore it holds always.
\end{proof}

\begin{remark}\label{5n11}
In order to have an absolutely continuous distribution of $
X_1,\dots,X_n $ for a given $ n $, the condition $ \Pr{ X_k = X_l } = 0
$ is needed only for $ k,l \in \{1,\dots,n\} $, $ k \ne l $.
\end{remark}

\section[]{\raggedright Main results}
\label{sect6}Recall Definitions \ref{4n2} (the independence condition) and \ref{2n2}
(the uniform distribution of a random countable set).

\begin{theorem}\label{6n1}
A random dense countable subset $ \{X_1,X_2,\dots\} $ of $ (0,1) $,
satisfying the independence condition, has the uniform distribution if
and only if
\begin{equation}\label{6n2}
\Pr{ B \cap \{X_1,X_2,\dots\} \ne \emptyset } = \begin{cases}
 0 &\text{if $ \mes(B)=0 $},\\
 1 &\text{if $ \mes(B)>0 $}
\end{cases}
\end{equation}
for all Borel sets $ B \subset (0,1) $. (Here `mes' is Lebesgue
measure.)
\end{theorem}

\begin{proof}
If it has the uniform distribution then we may assume that $
X_1,X_2,\dots $ are independent, uniform on $ (0,1) $, which makes
\eqref{6n2} evident.

Let \eqref{6n2} be satisfied. In order to prove that $ \{X_1,X_2,\dots\}
$ has the uniform distribution, it is sufficient to construct a
probability measure $ \mu $ on $ (0,1)^\infty \times (0,1)^\infty $ such
that the first marginal of $ \mu $ is the joint distribution of $
X_1,X_2,\dots $, the second marginal of $ \mu $ satisfies Conditions
(a), (b) of Main lemma \ref{3.2}, and $ \{x_1,x_2,\dots\} =
\{z_1,z_2,\dots\} $ for \almost{\mu} all $ \(
(x_1,x_2,\dots),(z_1,z_2,\dots) \) $.

To this end we construct recursively a consistent sequence of
probability measures $ \mu_n $ on $ (0,1)^\infty \times (0,1)^n $ (with
the prescribed first marginal) such that for all $ n $,
\begin{equation}\label{6n3}
\{ x_1,\dots,x_n \} \subset \{ z_1,\dots,z_{2n} \} \subset \{
x_1,x_2,\dots \}
\end{equation}
for \almost{\mu_{2n}} all $ \( (x_1,x_2,\dots),(z_1,\dots,z_{2n}) \) $, and
\begin{equation}\label{6n4}
z_{2n+1} \text{ is distributed uniformly and independent of }
z_1,\dots,z_{2n}
\end{equation}
w.r.t.\ $ \mu_{2n+1} $, and
\begin{equation}\label{6n45}
z_1, \dots, z_n \quad \text{are pairwise different}
\end{equation}
\almost{\mu_n} everywhere. This is sufficient since, first, \eqref{6n3}
implies $ \{x_1,x_2,\dots\} \linebreak[0]
= \{z_1,z_2,\dots\} $ for \almost{\mu} all $ \(
(x_1,x_2,\dots),(z_1,z_2,\dots) \) $ (here $ \mu $ is the measure
consistent with all $ \mu_n $); second, \ref{3.2}(a) is ensured by
\eqref{6n45}, Lemma \ref{5n10} and Remark \ref{5n11} (applied to $
(z_1,\dots,z_n,x_1,x_2,\dots) $ rather than $ (x_1,x_2,\dots) $); and
third, \eqref{6n4} implies \ref{3.2}(b).

We choose $ \mu_1 $ by means of Lemma \ref{3n7}.

For constructing $ \mu_{2n} $ we introduce Borel functions $ K_{2n} :
(0,1)^\infty \times (0,1)^{2n-1} \linebreak[0]
\to \{1,2,\dots\} $, $ Z_{2n} : (0,1)^\infty \times (0,1)^{2n-1} \to
(0,1) $ by
\begin{gather*}
K_{2n} (x,z) = \min \big\{ k : x_k \notin \{z_1,\dots,z_{2n-1}\} \big\}
 \, , \\
Z_{2n} (x,z) = x_{K_{2n} (x,z)} \, ;
\end{gather*}
of course, $ x=(x_1,x_2,\dots) $ and $ z=(z_1,\dots,z_{2n-1}) $.
Less formally, $ Z_{2n} (x,z) $ is the first of $ x_k $ different from
$ z_1,\dots,z_{2n-1} $. We define $ \mu_{2n} $ (consistent with $
\mu_{2n-1} $) such that
\[
z_{2n} = Z_{2n} \( (x_1,x_2,\dots), (z_1,\dots,z_{2n-1}) \)
\]
for \almost{\mu_{2n}} all $ (x,z) $. In other words, $ \mu_{2n} $ is the
distribution of
\[
\( (x_1,x_2,\dots), (z_1,\dots,z_{2n-1}, Z_{2n} (x,z) ) \)
\]
where $ (x,z) = \( (x_1,x_2,\dots), (z_1,\dots,z_{2n-1}) \) $ is
distributed $ \mu_{2n-1} $. Having \eqref{6n3} on the previous step,
$ \{x_1,\dots,x_{n-1}\} \subset \{ z_1,\dots,z_{2n-2} \} $, we conclude
that $ K_{2n} (x,z) \ge n $, thus, $ \{x_1,\dots,x_n\} \subset \{
z_1,\dots,z_{2n-1}, Z_{2n} (x,z) \} $, which ensures \eqref{6n3} on the
current step.

Finally, we choose $ \mu_{2n+1} $ by means of Prop.~\ref{5n9}.
\end{proof}

\begin{remark}
Assuming \eqref{3n8} instead of \eqref{6n2} we conclude that some part
of $ \{X_1,X_2,\dots\} $ (in the sense of \ref{2m5}) has the uniform
distribution. To this end we use only the `odd' part of the proof of
Th.~\ref{6n1}, that is, the construction of $ \mu_{2n+1} $.
\end{remark}

\begin{remark}
Assuming \eqref{4n3} instead of \eqref{6n2} we conclude that $
\{X_1,X_2,\dots\} $ is distributed as a part of a uniformly distributed
random set. To this end we use only the `even' part of the proof of
Th.~\ref{6n1}, that is, the construction of $ \mu_{2n} $, in combination
with Remark \ref{2m5}. 
\end{remark}

\begin{definition}\label{6n7}
A random dense countable subset $ \{X_1,X_2,\dots\} $ of $ (0,1) $,
satisfying the independence condition, is \emph{stationary,} if for
every $ a,b,c,d \in (0,1) $ such that $ b-a = d-c > 0 $, the two random
dense countable sets
\begin{gather*}
\big\{ x \in (0,1) : a+(b-a)x \in \{ X_1,X_2,\dots \} \big\}, \\
\big\{ x \in (0,1) : c+(d-c)x \in \{ X_1,X_2,\dots \} \big\}
\end{gather*}
are identically distributed.
\end{definition}

\begin{theorem}\label{6n8}
Every random dense countable subset  of $ (0,1) $, satisfying the
independence condition and stationary, has the uniform distribution.
\end{theorem}

\begin{proof}
First we prove \eqref{4n3}. Let $ B \subset (0,1) $ be a Borel set of
measure $ 0 $. Then for every $ x $, $ x \notin B+u $ for almost all $ u
$. Therefore for every $ \om $, $ \{ X_1(\om), X_2(\om), \dots \} \cap
(B+u) = \emptyset $ for almost all $ u $. By Fubini's theorem, $ \Pr{ \{
X_1, X_2, \dots \} \cap (B+u) = \emptyset } = 1 $ for almost all $ u
$. By stationarity, this probability does not depend on $ u $ as long as
$ B \subset (0,\frac12) $ and $ u \in [0,\frac12] $, or $ B \subset
(\frac12,1) $ and $ u \in [-\frac12,0] $; \eqref{4n3} follows.

Second, Prop.~\ref{4n4} gives us a function $ r : (0,1) \to [0,\infty]
$. By Remark \ref{4n5} and stationarity, this function (or rather, its
equivalence class) is shift invariant; thus, $ r(x)=\infty $ for almost
all $ x $, which implies \eqref{3n8}. It remains to apply
Th.~\ref{6n1}.
\end{proof}

It is well-known (see \cite[2.9.12]{KS}) that for almost every Brownian
path $ w : [0,\infty) \to \R $ each local minimizer (that is, $ x \in
(0,\infty) $ such that $ w(y) \ge w(x) $ for all $ y $ close enough to $
x $) is a strict local minimizer (it means, $ w(y) > w(x) $ for all $
y $ close enough to $ x $, except for $ x $ itself), and all local
minimizers are a dense countable set.

\begin{lemma}\label{6n9}
There exist Borel functions $ X_1, X_2, \dots : C[0,1] \to (0,1) $ such
that for almost every Brownian path $ w $ (that is, for almost all $ w
\in C[0,1] $ w.r.t.\ the Wiener measure) the set $ \{ X_1(w), X_2(w),
\dots \} $ is equal to the set of all local minimizers of $ w $.
\end{lemma}

I give two proofs.

\begin{proof}[First proof]
We take a sequence of intervals $ (a_k,b_k) \subset (0,1) $ that are a
base of the topology, and define $ X_k(w) $ as the (global) minimizer of
$ w $ on $ [a_k,b_k] $ whenever it is unique.
\end{proof}

\begin{proof}[Second proof]
It is observed by Kendall \cite[Th.~3.4]{Ke00} that well-known
selection theorems (see \cite[Th.~(18.10)]{Ke95}) can be used for
constructing $ X_1,X_2,\dots $ provided that the set of pairs
\[
\{ (w,x) : x \text{ is a local minimizer of } w \}
\]
is a Borel subset of $ C[0,1] \times (0,1) $. It remains to note that
for every $ \eps>0 $ the set of pairs $ (w,x) $ such that $ w(x) =
\min_{[x-\eps,x+\eps]} w $ is closed.
\end{proof}

\begin{theorem}\label{6n10}
The random dense countable set of all local minimizers of a Brownian
motion on $ (0,1) $ has the uniform distribution.
\end{theorem}

\begin{proof}
We start with Lemma \ref{6n9}, note that the independence condition and
stationarity hold, and apply Theorem \ref{6n8}.
\end{proof}

All the arguments can be generalized to higher dimensions and applied to
the other, percolation-related, models mentioned in Introduction; see
\ref{9n6}.

\section[]{\raggedright Borelogy, the new framework}
\label{sect7}First of all, two quotations.

\begin{quote}
It has long been recognized in diverse areas of mathematics that in many
important cases such quotient spaces $ X/E $ cannot be viewed as
reasonable subsets of Polish spaces and therefore the usual methods of
topology, geometry, measure theory, etc., are not directly applicable
for their study. Thus they are often referred to as \emph{singular
spaces.}\hfill\mbox{Kechris \cite[\S 2]{Ke99}.}
\end{quote}

\begin{quote}
In this theory, a differential structure of some set $ X $ is defined as
the set of all the ``differentiable parametrizations'' of $ X $ [\dots]
The set of these chosen parametrizations is called a \emph{diffeology}
of $ X $, and its elements are called the \emph{plots} of the
diffeology. [\dots] In other words, a diffeology of $ X $ says how to
``walk'' differentiably into $ X $.\hfill\mbox{Iglesias-Zemmour
\cite[beginning of Chapter~1]{IZ}.}
\end{quote}

I propose \emph{borelogy}, a solution of the following `ideological
equation',
\[
\frac{ \text{borelogical space} }{ \text{standard Borel space} } =
\frac{ \text{diffeological space} }{ \text{differentiable manifold} } \,
.
\]
The set $ \DCS(0,1) $ of all dense countable subsets of the interval $
(0,1) $ is an example of a singular space in the sense of Kechris (see $
F_2 $ in \cite[\S 8]{Ke99}). It is also an example of a borelogical
space, see \ref{7n7} below. By the way, some singular spaces are
diffeological spaces, see \cite[1.15]{IZ}.

\begin{definition}\label{1.1}
A \emph{borelogy} on a set $ V $ is a set $ \BB $ of maps $ \R \to V $,
such that the following three conditions are satisfied:

(a) for every $ b \in \BB $,
\begin{equation}\label{1.2}
\text{the set } \{ (x,y) : b(x) = b(y) \} \text{ is a Borel subset of $
  \R^2 $} \, ;
\end{equation}

(b) for every Borel function $ f : \R \to \R $ and every $ b \in \BB $
their composition $ b(f(\cdot)) $ belongs to $ \BB $;

(c) there exists $ b \in \BB $ such that $ b(\R) = V $ and the set of
compositions $ b(f(\cdot)) $ (where $ f $ runs over all Borel functions
$ \R \to \R $) is the whole $ \BB $.
\end{definition}

\begin{definition}
A \emph{borelogical space} is a pair $ (V,\BB) $ of a set $ V $ and a
borelogy $ \BB $ on $ V $. Elements of $ \BB $ are called \emph{plots}
of the borelogical space.
\end{definition}

\begin{lemma}\label{1.4}
Let $ b : \R \to V $ be a surjective map satisfying \eqref{1.2}, and $ \BB
$ be the set of compositions $ b(f(\cdot)) $ for all Borel $ f : \R \to
\R $. Then $ \BB $ is a borelogy on $ V $.
\end{lemma}

\begin{proof}
(a): The set $ \{ (x,y) : b(f(x))=b(f(y)) \} $ is the inverse image of
the Borel set $ \{ (x,y) : b(x)=b(y) \} $ under the Borel map $ (x,y)
\mapsto (f(x),f(y)) $, therefore, a Borel set.

(b): $ b(f(g(\cdot))) $ belongs to $ \BB $, since $ f(g(\cdot)) $ is a
Borel function.

(c): The given $ b $ fits by construction.
\end{proof}

We see that every surjective map $ b : \R \to V $ satisfying \eqref{1.2}
generates a borelogy. If $ f : \R \to \R $ is a Borel isomorphism (that
is, $ f $ is invertible, and $ f $, $ f^{-1} $ both are Borel functions)
then $ b $ and $ b(f(\cdot)) $ generate the same borelogy.

Recall that an uncountable standard Borel space may be defined as a
measurable space, Borel isomorphic to $ \R $. Every uncountable standard
Borel space $ S $ may be used instead of $ \R $ in the definition of a
borelogy. Such \based{S} borelogies are in a natural one-to-one
correspondence with \based{\R} borelogies. The correspondence is
established via an isomorphism between $ S $ and $ \R $, but does not
depend on the choice of the isomorphism.

According to \ref{1.1}(c), every borelogy $ \BB $ contains a \emph{generating
plot,} that is, $ \BB $ is generated by some $ b : \R \to V $. Such $ b $
establishes a bijective correspondence between $ V $ and the quotient
set $ \R / E_b $, where $ E_b = \{ (x,y) : b(x)=b(y) \} \subset \R^2 $
is the relevant equivalence relation (a Borel equivalence relation, due
to \eqref{1.2}). Every $ b' \in \BB $ is of the
form $ b'(\cdot) = b(f(\cdot)) $ for some Borel $ f : \R \to \R
$. Clearly, $ (x,y) \in E_{b'} $ if and only if $ (f(x),f(y)) \in E_b
$. In terms of \cite[\S 3]{Ke99} it means that $ f $ is a Borel
reduction of $ E_{b'} $ into $ E_b $, and $ \R / E_{b'} $ has Borel
cardinality at most that of $ \R / E_b $. It may happen that also $ b' $
generates $ \BB $. Then one says that $ \R / E_{b'} $ and $ \R / E_b $
have the same Borel cardinality \cite[\S 3]{Ke99}. We see that every
borelogical space has its (well-defined) Borel cardinality.

Recall also that a standard Borel space is either an uncountable
standard Borel space (discussed above), or a finite or countable set
equipped with the \sif\ of all subsets. Let $ (V,\B) $ be a standard
Borel space (here $ \B $ is a given \sif\ of subsets of $ V $). We turn
it into a borelogical space $ (V,\BB) $ where $ \BB $ consists of all
Borel maps $ b : \R \to V $. Such a borelogical space will be called
\emph{nonsingular.} This way, standard Borel spaces may be treated as a
special case of borelogical spaces. In terms of the equivalence relation
$ E_b $ corresponding to a generating plot $ b $, the borelogical space
is nonsingular if and only if $ E_b $ is smooth (or tame), see \cite[\S
6]{Ke99}. Otherwise, the borelogical space will be called
\emph{singular}.

A nonempty finite or countable set $ V $ carries one and only one
borelogy $ \BB $, and $ (V,\BB) $ is nonsingular.

Pinciroli \cite[Def.~1.2]{Pi} defines a quotient Borel space as a couple
$ (S,E) $ where $ S $ is a standard Borel space and $ E $ is a Borel
equivalence relation on $ S $ whose equivalence classes are
countable. He stipulates that the underlying set of $ (S,E) $ is the
quotient set $ S/E $. Clearly, every quotient Borel space is a
borelogical space. On the other hand, every borelogical space is of the
form $ S/E $, however, a freedom is left in the choice of $ S $ and $ E
$ (see also Example \ref{7n8}), and $ E $ need not have countable
equivalence classes.

\begin{example}\label{1.5}
The set $ V = \R / \Q $ (reals modulo rationals) consists of equivalence
classes $ \Q + x = \{ q+x : q\in\Q \} $ for all $ x \in \R $. The
natural map $ \R \to V $, $ x \mapsto \Q+x $, generates a borelogy $ \BB
$ on $ V $ (by Lemma \ref{1.4}), and $ (V,\BB) $ is singular (see Remark
\ref{2.R}).
\end{example}

\begin{example}\label{1.6}
The set $ V = \CS(0,1) $ of all countable subsets of the interval $
(0,1) $ is the image of the set $ S = (0,1)^\infty_{\ne} \subset
(0,1)^\infty $ of all  sequences $ (s_1,s_2,\dots) $ of pairwise
different numbers $ s_k \in (0,1) $ under the natural map $ b : S \to V
$, $ b(s_1,s_2,\dots) = \{s_1,s_2,\dots\} $.

The set $ S $ is a Borel subset of $ (0,1)^\infty $ (since $ \{
(s_1,s_2,\dots) : s_k \ne s_l \} $ is open whenever $ k \ne l $),
therefore, a standard Borel space, see \cite[Sect.~12.B]{Ke95}.

The map $ b : S \to V $ satisfies \eqref{1.2}, that is, the set of all
pairs $ (s,s') = ((s_1,s_2,\dots),(s'_1,s'_2,\dots)) $ such that $ b(s)
= b(s') $ is a Borel subset of $ S \times S $, since
\[
b(s) = b(s') \equi \forall k \> \exists l \; (s_k = s'_l) \;\>\&\;\>
\forall k \> \exists l \; (s'_k = s_l) \, .
\]
By an evident generalization of Lemma \ref{1.4}, $ b $ generates on $ V
$ an \based{S} borelogy, which turns $ V $ into a borelogical space. It
is singular (see Remark \ref{2.R}).
\end{example}

\begin{example}\label{7n7}
The borelogical space $ \DCS(0,1) $ of all \emph{dense} countable
subsets of $ (0,1) $ is defined similarly. It is singular (see Remark
\ref{2.R}).
\end{example}

\begin{example}\label{7n8}
The borelogical space $ \FCS(0,1) $ of all \emph{finite or countable}
subsets of $ (0,1) $ is defined similarly. We may adopt finite sets by
replacing $ (0,1)^\infty_{\ne} $ with $ \cup_{n=0,1,2,\dots;\infty}
(0,1)^n_{\ne} $, where $ (0,1)^0_{\ne} $ contains a single element (the
empty sequence) whose image in $ \FCS(0,1) $ is the empty
set. Alternatively we may adopt finite sets by replacing $
(0,1)^\infty_{\ne} $ with $ (0,1)^\infty $, thus permitting equal
numbers in the sequences; however, in this case we should bother about
the empty set as an element of $ \FCS(0,1) $.

We may also consider the set $ S $ of all discrete finite positive Borel
measures on $ (0,1) $ (`discrete' means existence of a countable set of
full measure) together with the equivalence relation $ E $ of mutual
absolute continuity. Once again, $ S/E = \FCS(0,1) $. \emph{Sketch of
the proof.} On one hand, a sequence $ (s_1,s_2,\dots) $ leads to a
discrete measure $ A \mapsto \sum_{k:s_k\in A} 2^{-k} $. On the
other hand, a discrete measure $ \mu $ leads to the sequence of its
atoms, the most massive atom being the first and so on. (If several
atoms are equally massive, the leftmost one is the first.)
\end{example}

As usual, we often say `a borelogical space $ V $' rather than `a
borelogical space $ (V,\BB) $'.

\begin{definition}
Let $ V,W $ be borelogical spaces.

(a) A \emph{morphism} of $ V $ to $ W $ is a map $ f : V \to W $ such
that for every plot $ b $ of $ V $ the map $ f(b(\cdot)) $ is a plot of
$ W $.

(b) An \emph{isomorphism} between $ V $ and $ W $ is an invertible map $
f : V \to W $ such that $ f $ and $ f^{-1} $ are morphisms.
\end{definition}

Choosing generating plots $ b $ for $ V $ and $ b' $ for $ W $ we
observe that $ f : V \to W $ is a morphism if and only if $ f(b(\cdot))
= b'(g(\cdot)) $ for some Borel $ g : \R \to \R $.
\[
\xymatrix{
 \R \ar[r]^{g} \ar[d]_{b} & \R \ar[d]^{b'}
\\
 V \ar[r]_{f} & W
}
\]
(Compare it with a Borel morphism of Borel equivalence relations
\cite[p.~1]{Pi}.)
If a morphism $ f : V \to W $ is injective (that is, $ x_1 \ne x_2 $
implies $ f(x_1) \ne f(x_2) $) then
\[
(x_1,x_2) \in E_b \equi (g(x_1),g(x_2)) \in E_{b'} \, ,
\]
which means that $ g $ is a Borel reduction of $ E_b $ into $ E_{b'} $,
and the Borel cardinality of $ V $ is at most that of $ W $. It follows
that isomorphic borelogical spaces have the same Borel cardinality. (Is
the converse true? I do not know.) Existence of a continuum of
(different) Borel cardinalities, mentioned in \cite[\S 6]{Ke99}, implies
existence of a continuum of mutually nonisomorphic borelogical
spaces. Of course, they are singular; all nonsingular borelogical spaces
of the same cardinality (finite, countable or continuum) are isomorphic.

The Borel cardinality of $ \R / \Q $ is well-known as $ E_0 $ \cite[\S
3]{Ke99}. The Borel cardinality of $ \CS(0,1) $ is well-known as $ F_2
$ \cite[\S 8]{Ke99}.

\begin{definition}\label{1.8}
The \emph{product} of two borelogical spaces $ (V_1,\BB_1), (V_2,\BB_2)
$ is the borelogical space $ (V_1\times V_2, \BB_1\times\BB_2) $ where $
\BB_1\times\BB_2 $ consists of all maps $ \R \to V_1 \times V_2 $ of the
form $ b_1 \times b_2 $, that is, $ x \mapsto \( b_1(x),b_2(x) \) $,
where $ b_1 \in \BB_1 $, $ b_2 \in \BB_2 $.
\end{definition}

\begin{lemma}
Definition \ref{1.8} is correct, that is, $ \BB_1\times\BB_2 $ is a
borelogy on $ V_1 \times V_2 $.
\end{lemma}

\begin{proof}
We check the three conditions of \ref{1.1}.

(a) the set $ \{ (x,y) : \( b_1(x),b_2(x) \) = \( b_1(y),b_2(y) \) \} =
\{ (x,y) : b_1(x)=b_1(y) \} \cap \{ (x,y) : b_2(x)=b_2(y) \} $ is the
intersection of two Borel sets;

(b) $ (b_1\times b_2) (f(\cdot)) = \( b_1(f(\cdot)), b_2(f(\cdot)) \) =
b_1(f(\cdot)) \times b_2(f(\cdot)) \in \BB_1 \times \BB_2 $;

(c) choosing generating plots $ b_1 : \R \to V_1 $, $ b_2 : \R \to V_2 $
and a Borel isomorphism $ g : \R \to \R^2 $, $ g(\cdot) = \( g_1(\cdot),
g_2(\cdot) \) $, we define $ b : \R \to V_1 \times V_2 $ by $ b(x) = \(
b_1(g_1(x)), b_2(g_2(x)) \) $ and note that every element of $ \BB_1
\times \BB_2 $ is of the form $ x \mapsto \( b_1(f_1(x)), b_2(f_2(x)) \)
= \( b_1(g_1(f(x))), b_2(g_2(f(x))) \) = b(f(x)) $ where $ f(x) = g^{-1}
\( ( f_1(x), f_2(x) ) \) $.
\end{proof}

\begin{remark}
Having usual (\based{\R}) generating plots $ b_1 : \R \to V_1 $, $ b_2
: \R \to V_2 $ we get immediately an \based{\R^2} generating plot $ b :
\R^2 \to V_1 \times V_2 $, $ b(x,y) = \( b_1(x), b_2(y) \) $. However,
in order to get a generating plot $ \R \to V_1 \times V_2 $ we need a
Borel isomorphism between $ \R $ and $ \R^2 $.

More generally, having generating plots $ b_1 : S_1 \to V_1 $, $ b_2 :
S_2 \to V_2 $ we get immediately a generating plot $ b : S_1 \times S_2
\to V_1 \times V_2 $, $ b(x,y) = \( b_1(x), b_2(y) \) $; here $ S_1, S_2
$ are standard Borel spaces.
\end{remark}

Compare Def.~\ref{1.8} with \cite[Def.~1.3]{Pi}.

\begin{example}
Defining a borelogical space $ \R^2 / \Q^2 $ similarly to \ref{1.5} we
get (up to a natural isomorphism)
\[
(\R/\Q) \times (\R/\Q) = \R^2/\Q^2 \, .
\]
\end{example}

\begin{example}
Defining a borelogical space $ \CS[a,b) $ for any $ [a,b) \subset \R $
similarly to \ref{1.6} we get (up to a natural isomorphism)
\[
\CS[a,b) \times \CS[b,c) = \CS[a,c)
\]
and the same for $ \DCS $ and $ \FCS $.
\end{example}

\section[]{\raggedright Probability measures on singular spaces}
\label{sect8}Each borelogical space $ V $ carries a \sif\ $ \Si $ consisting of all $
A \subset V $ such that $ b^{-1} (A) $ is a Borel subset of $ \R $ for
every plot $ b : \R \to V $. Choosing a generating plot $ b $ we see
that $ A \in \Si $ if and only if $ b^{-1} (A) $ is a Borel set. (That
is, $ \Si $ is the quotient \sif\ \cite[Sect.~5.1]{Sri}, see also
\cite[Sect.~1]{Pi}.) If $ V $ is nonsingular then $ \Si $ is its given
Borel \sif\ (since the identical map $ V \to V $ is a plot), thus, $
(V,\Si) $ is a standard Borel space. If $ V $ is singular then it may
happen that $ (V,\Si) $ is still a standard Borel space (which can be
shown by means of a well-known counterexample \cite[5.1.7]{Sri}), but in
such cases as $ \R / \Q $ and $ \DCS(0,1) $ it is not.

In my opinion, a notion defined via a singular space can be useful in
probability theory only if it admits an equivalent definition in terms
of standard spaces. A quote from Pinciroli \cite[p.~2]{Pi}: ``[\dots]
the `right' notion of Borelness for [\dots] functions between quotient
Borel spaces is \emph{not} the usual one from the context of measurable
spaces and maps: here again we want to exploit the original standard
Borel structures.''

Three examples follow.

First, it may be tempting to define a random element of a borelogical
space $ V $ as a measurable map from a standard probability space $ \Om
$ to the (nonstandard) measurable space $ (V,\Si) $. However, I prefer
to define a random element of $ V $ as a map $ \Om \to V $  of the form $
b(X(\cdot)) $ where $ X : \Om \to \R $ is a (usual) random variable, and
$ b : \R \to V $ is a plot. Are these two definitions equivalent? I do
not know. Every $ b(X(\cdot)) $ is \measurable{\Si}, but I doubt that
every \measurable{\Si} map is of the form $ b(X(\cdot)) $. For Borel
maps the answer is negative, but for equivalence classes it may be
different.

Second, it may be tempting to define the distribution of a random
element $ b(X(\cdot)) $ as the corresponding probability measure on the
(nonstandard) measurable space $ (V,\Si) $. Then two random elements may
be treated as identically distributed if their distributions are
equal. An equivalent definition in terms of standard spaces will be
given (Th.~\ref{2.2}, Def.~\ref{2.3}).

Third, it may be tempting to say that two random elements $ b(X) $ and $
b(Y) $ are independent if $ \Pr{ b(X) \in A, \, b(Y) \in B } = \Pr{ b(X)
\in A } \Pr{ b(Y) \in B } $ for all $ A,B \in \Si $. However, I prefer a
different, nonequivalent definition (see Def.~\ref{2.16} and
Counterexample \ref{2.17}).

Recall that a standard probability space (known also as a
Lebesgue-Rokhlin space) is a probability space isomorphic $ (\modO) $ to
an interval with the Lebesgue measure, a finite or countable collection
of atoms, or a combination of both. Every probability measure on a
standard Borel space turns it (after completion, that is, adding all
negligible sets to the \sif) into a standard probability space, see
\cite[Sect.~17.F]{Ke95} or \cite[Th.~3.4.23]{Sri}.

\begin{definition}\label{8n1}
Let $ V $ be a borelogical space and $ \Om $ a standard probability
space. A \emph{\valued{V\!} random variable} on $ \Om $ (called also a
random element of $ V $) is an equivalence class of maps $ X : \Om \to V
$ representable in the form $ X(\cdot) = b(Y(\cdot)) $ for some plot $ b
: \R \to V $ and some (usual) random variable $ Y : \Om \to \R
$. (Equivalence means equality almost everywhere on $ \Om $.)
\end{definition}

Choosing a generating plot $ b : \R \to V $ we observe that every random
element of $ V $ is of the form $ b(Y(\cdot)) $ with the chosen $ b $
and arbitrary $ Y $, since $ b(f(Y(\cdot))) $ is of this form.

\begin{theorem}\label{2.2}
Let $ (\Om_1,\F_1,P_1), (\Om_2,\F_2,P_2) $ be standard probability
spaces, $ V $ a borelogical space, and $ X_1, X_2 $ be \valued{V} random
variables on $ \Om_1, \Om_2 $ respectively. Then the following two
conditions are equivalent.

(a) $ \Pr{ X_1 \in A } = \Pr{ X_2 \in A } $ for all $ A \in \Si $;

(b) there exists a probability measure $ P $ on $ \Om_1 \times \Om_2 $
whose marginals are $ P_1, P_2 $, such that $ X_1(\om_1) = X_2(\om_2) $
for \almost{P} all $ (\om_1,\om_2) $.
\end{theorem}

The proof is given after \ref{2.9}.

Theorem \ref{2.2} shows that items (b1), (b2) of the following
definition are basically the same.

\begin{definition}\label{2.3}
Let $ V $ be a borelogical space.

(a) Two \valued{V} random variables are \emph{identically distributed,}
if they satisfy the equivalent conditions (a), (b) of Theorem \ref{2.2}.

(b1) A \emph{distribution} on $ V $ is an equivalence class of
\valued{V} random variables on the probability space $ \Om = (0,1) $
(with Lebesgue measure); here random variables are treated as equivalent
if they are identically distributed.

(b2) A \emph{distribution} on $ V $ is a probability measure on the
(generally, nonstandard) measurable space $ (V,\Si) $, representable in
the form $ \Pr{X\in\cdot} $ for some \valued{V} random variable $ X $.

(c) A distribution on $ V $ is called an \OI\ distribution, if
it ascribes to all sets of $ \Si $ the probabilities $ 0,1 $ only.
\end{definition}

\begin{remark}\label{2.31}
If $ X(\cdot) = b(Y(\cdot)) $ has a \OI\ distribution and $ Y' $ has a
distribution absolutely continuous w.r.t.\ the distribution of $ Y $,
then $ X' = b(Y') $ and $ X $ are identically distributed.
\end{remark}

\begin{example}\label{2.32}
Continuing Example \ref{1.5} we consider the borelogical space $ \R / \Q
$ and its generating plot $ b : \R \to \R/\Q $, $ b(x) = \Q + x $. Every
random variable $ Y : \Om \to \R $ leads to a \valued{\R/\Q}
random variable $ X = \Q + Y $. If the distribution of $ Y $ is
absolutely continuous then $ X $ has a \OI\ distribution (since for
every \invariant{\Q} Borel $ B \subset \R $ it is well-known that either
$ B $ or $ \R \setminus B $ is of Lebesgue measure $ 0 $). All
absolutely continuous distributions on $ \R $ correspond to a single
distribution on $ \R/\Q $. This special \OI\ distribution on $ \R/\Q $
may be called \emph{the uniform distribution on $ \R/\Q $.}
\end{example}

\begin{corollary}\label{2.33}
For every two absolutely continuous probability measures $ \mu,\nu $ on
$ \R $ there exist random variables $ X,Y $ distributed $ \mu,\nu $
respectively and such that the difference $ X-Y $ is a.s.\ a (random)
rational number.
\end{corollary}

If $ X_1(\cdot) = b(Y_1(\cdot)) $, $ X_2(\cdot) = b(Y_2(\cdot)) $ and $
Y_1,Y_2 $ are identically distributed then, of course, $ X_1,X_2 $ are
identically distributed. The converse does not hold (without
an appropriate enlargement of probability spaces), see below.

\begin{counterexample}
There exist two identically distributed \valued{\R/\Q} random variables
$ X_1, X_2 : \Om \to \R / \Q $ that are \emph{not} of the form $ X_1 =
\Q + Y_1 $, $ X_2 = \Q + Y_2 $ where $ Y_1,Y_2 : \Om \to \R $ are
identically distributed.
\end{counterexample}

\begin{proof}
(See also \ref{3n9}.) We take $ \Om = (0,1) $ with Lebesgue measure and
define for $ \om \in \Om $
\[
X_1(\om) = \Q + \om \, , \quad X_2(\om) = \Q + \sqrt2 \, \om \, .
\]
Let $ Y_1 : \Om \to \R $ satisfy $ X_1 = \Q + Y_1 $ a.s.; I claim that
the distribution of $ Y_1 $ necessarily has an integer-valued
density. Similarly, I claim that every $ Y_2 : \Om \to \R $ satisfying $
X_2 = \Q + Y_2 $ a.s.\ necessarily has a density that takes on the
values $ 0, 1/\sqrt2,  2/\sqrt2,  3/\sqrt2, \dots $ only. Clearly, such
$ Y_1, Y_2 $ cannot be identically distributed; it remain to prove the
first claim (the second claim is similar).

We have $ \Q + \om = \Q + Y_1(\om) $, that is, $ Y_1(\om) - \om \in \Q
$. We partition $ (0,1) $ into countably many measurable sets $ A_q = \{
\om : Y_1(\om) - \om = q \} $ for $ q \in \Q $ and observe that $ \Pr{
Y_1 \in B } = \sum_{q\in\Q} \mes \{ \om \in A_q : \om + q \in B \} =
\int_B f(x) \, \D x $ where $ f(x) $ is the number of $ q \in \Q $ such
that $ x-q \in A_q $.
\end{proof}

\begin{example}\label{2.37}
Continuing Example \ref{1.6} we consider the borelogical space $
\CS(0,1) $ and its generating plot $ b : (0,1)^\infty_{\ne} \to \CS(0,1)
$. An \valued{(0,1)^\infty_{\ne}} random variable $ Y $ is nothing but a
sequence of random variables $ Y_1,Y_2,\dots : \Om \to (0,1) $ such that
$ \Pr{ Y_k = Y_l } = 0 $ for $ k \ne l $. Every such $ Y $ leads to a
\valued{\CS(0,1)} random variable $ X = b(Y) $, that is, a random
countable set $ X(\om) = \{ Y_1(\om), Y_2(\om), \dots \} $; see also
\eqref{2n0}. If $ Y_k $ are independent then $ X $ has a \OI\
distribution by the Hewitt-Savage \OI\ law. Using \ref{2.31}, similarly
to \ref{2.33}, if $ Y_k $ are independent and the distribution of $
Y'=(Y'_1,Y'_2,\dots) $ is absolutely continuous w.r.t.\ the distribution
of $ Y=(Y_1,Y_2,\dots) $, then there exists a joining between $ Y $ and
$ Y' $ such that $ \{ Y_1(\om), Y_2(\om), \dots \} = \{ Y'_1(\om),
Y'_2(\om), \dots \} $ a.s.

According to Main lemma \ref{3.2}, a wide class of probability
distributions on $ (0,1)^\infty_{\ne} $ (many of them being mutually
singular) corresponds to a single \OI\ distribution on $ \CS(0,1) $,
called uniform according to \ref{2n2}. (See also \ref{2.32}.)
\end{example}

\begin{remark}\label{2.R}
Every \OI\ distribution on a standard Borel space (that is, a
nonsingular borelogical space) is concentrated at a single
point. Therefore, existence of a \OI\ distribution that does not charge
points implies singularity of a borelogical space. (See also
\cite[Remark~3.3]{Pi}.) In particular, $ \R/\Q $ and $ \DCS(0,1) $ are
singular (recall \ref{2.32} and \ref{2.37}).
\end{remark}

The proof of Theorem \ref{2.2} is based on Kellerer's Theorem \ref{2n4};
recall it: $ S_{\mu_1,\mu_2} (B) = I_{\mu_1,\mu_2} (B) $. The theorem
holds for all standard Borel spaces $ \X_1, \X_2 $ and Borel sets $ B
\subset \X_1 \times \X_2 $; we apply it to $ \X_1 = \X_2 = \R $ and
specialize the Borel set as follows.

Let a Borel set $ E \subset \R^2 $ be (the graph of) an equivalence
relation on $ \R $; we introduce the \sif\ $ \Ec $ of all Borel sets $ A
\subset \R $ that are saturated (invariant) in the sense that
\[
(x,y) \in E \imply ( x \in A \equi y \in A ) \qquad \text{for } x,y \in
\R \, .
\]

\begin{lemma}\label{2.7}
For all probability measures $ \mu_1,\mu_2 $ on $ \R $,
\[
I_{\mu_1,\mu_2} (E) = 1 - \sup_{A\in\Ec} | \mu_1(A) - \mu_2(A) | \, .
\]
\end{lemma}

\begin{proof}
First, ``$\le$'': we have $ E \subset (A\times\R) \cup (\R\times
\overline A) $ for $ A \in \Ec $ (here $ \overline A = \R \setminus A
$), therefore $ I_{\mu_1,\mu_2} (E) \le \mu_1(A) + \mu_2(\overline A) =
1 - \( \mu_2(A)-\mu_1(A) \) $. Similarly, $ I_{\mu_1,\mu_2} (E) \le 1 -
\( \mu_1(A)-\mu_2(A) \) $.

Second, ``$\ge$''. Let Borel $ B_1,B_2 $ satisfy $ E \subset
(B_1\times\R) \cup (\R\times B_2) $, then $ \overline B_1 \times
\overline B_2 $ does not intersect $ E $. It follows (see
\cite[Th.~4.4.5]{Sri} or \cite[Exercise (14.14)]{Ke95}) that there
exists $ A \in \Ec $ such that $ \overline B_1 \subset A $, $ \overline
B_2 \subset \overline A $. We have $ \mu_1(B_1) + \mu_2(B_2) \ge
\mu_1(\overline A) + \mu_2(A) = 1 - \( \mu_1(A)-\mu_2(A) \) \ge 1 -
\sup_{A\in\Ec} | \mu_1(A)-\mu_2(A) | $.
\end{proof}

\begin{remark}\label{2.8}
If $ \mu_1,\mu_2 $ are finite positive (not just probability) measures
such that $ \mu_1(\R) = \mu_2(\R) $ then $ I_{\mu_1,\mu_2} (E) =
\mu_1(\R) - \sup_{A\in\Ec} | \mu_1(A) - \mu_2(A) | $.
\end{remark}

\begin{lemma}\label{2.9}
The supremum $ S_{\mu_1,\mu_2} (E) $ is reached for all probability
measures $ \mu_1,\mu_2 $ on $ \R $.
\end{lemma}

\begin{proof}
We check the condition of Lemma \ref{2n6}. Let $ \nu $ be a positive
measure on $ E $ with marginals $ \nu_1 \le \mu_1 $, $ \nu_2 \le \mu_2
$. By Theorem \ref{2n4} it is sufficient to prove that $
I_{\mu_1-\nu_1,\mu_2-\nu_2} (E) = I_{\mu_1,\mu_2} (E) - \nu(E) $. By
\ref{2.7} and \ref{2.8} it boils down to the equality
\[
\sup_{A\in\Ec} | (\mu_1-\nu_1) (A) - (\mu_2-\nu_2) (A) | =
\sup_{A\in\Ec} | \mu_1(A) - \mu_2(A) | \, .
\]
It remains to note that $ \nu_1(A) = \nu_2(A) $ for all $ A \in \Ec $.
\end{proof}

\begin{proof}[Proof of Theorem \textup{\ref{2.2}}]
(b) \imp (a):
$ P_1 \( \{ \om_1 : X_1(\om_1) \in A \} \) = P \( \{ (\om_1,\om_2) :
X_1(\om_1) \in A \} \) = P \( \{ (\om_1,\om_2) : X_2(\om_2) \in A \} \)
= P_2 \( \{ \om_2 : X_2(\om_2) \in A \} \) $.

(a) \imp (b):
We choose a generating plot $ b : \R \to V $ and random variables $ Y_1
: \Om_1 \to \R $, $ Y_2 : \Om_2 \to \R $ such that $ X_1(\cdot) =
b(Y_1(\cdot)) $, $ X_2(\cdot) = b(Y_2(\cdot)) $. The Borel set $ E = \{
(x,y) : b(x) = b(y) \} \subset \R^2 $ is an equivalence relation. Let $
A \in \Ec $ (that is, $ A $ is a saturated Borel set), then $ A = b^{-1}
\( b(A) \) $ and $ b(A) \in \Si $. It is given that $ P_1 \( X_1 \in
b(A) \) = P_2 \( X_2 \in b(A) \) $, that is, $ P_1 \( Y_1 \in A \) =P_2
\( Y_2 \in A \) $. Denoting by $ \mu_1,\mu_2 $ the distributions of $
Y_1, Y_2 $ respectively, we see that $ \mu_1(A) = \mu_2(A) $ for all $ A
\in \Ec $.

By Lemma \ref{2.7}, $ I_{\mu_1,\mu_2} (E) = 1 $. By Theorem \ref{2n4}, $
S_{\mu_1,\mu_2} (E) = 1 $. Lemma \ref{2.9} gives us a probability
measure $ \mu $ on $ \R^2 $ with the marginals $ \mu_1,\mu_2 $ such that
$ \mu(E) = 1 $.

We consider the conditional distribution $ P_{1,x} $ of $ \om_1 \in
\Om_1 $ given $ Y_1(\om_1) = x $ (its existence is ensured by
standardness of $ \Om_1 $); $ P_{1,x} $ is a probability measure on $
\Om_1 $ for \almost{\mu_1} all $ x \in \R $, and $ \int_{\R} P_{1,x} \,
\mu_1(\D x) = P_1 $. The same holds for $ P_{2,y} $. We construct a
probability measure $ P $ on $ \Om_1 \times \Om_2 $ by
\[
P = \int_{\R^2} ( P_{1,x} \times P_{2,y} ) \, \mu(\D x \D y) \, .
\]
The first marginal of $ P $ is $ \int P_{1,x} \, \mu(\D x \D y) = \int
P_{1,x} \, \mu_1(\D x) = P_1 $; the second marginal of $ P $ is $ P_2
$. Also, $ P \( \{ (\om_1,\om_2) : ( Y_1(\om_1),Y_2(\om_2) ) \in E \} \)
= \mu(E) = 1 $. Thus, for \almost{P} all $ (\om_1,\om_2) $ we have $ (
Y_1(\om_1),Y_2(\om_2) ) \in E $, therefore, $ b(Y_1(\om_1)) =
b(Y_2(\om_2)) $, that is, $ X_1(\om_1) = X_2(\om_2) $.
\end{proof}

\begin{definition}\label{2.16}
Let $ V $ be a borelogical space and $ \Om $ a standard probability
space. Two \valued{V} random variables $ X_1,X_2 : \Om \to V $ are
\emph{independent,} if there exist independent random variables $
Y_1,Y_2 : \Om \to \R $ and a plot $ b : \R \to V $ such that $
X_1(\cdot) = b(Y_1(\cdot)) $ and $ X_2(\cdot) = b(Y_2(\cdot)) $.
Independence of three or more \valued{V} random variables is defined
similarly, as well as independence of random elements of different
borelogical spaces.
\end{definition}

A given generating plot $ b : \R \to V $ can be used always. Also, the
case $ X_1(\cdot) = b_1(Y_1(\cdot)) $, $ X_2(\cdot) = b_2(Y_2(\cdot)) $
reduces to a single $ b $.

If $ X_1,X_2 : \Om \to V $ are independent then $ \Pr{ X_1 \in A_1, X_2
\in A_2 } = \Pr{ X_1 \in A_1 } \Pr{ X_2 \in A_2 } $ for all $ A_1,A_2
\in \Si $. The converse is generally wrong.

\begin{counterexample}\label{2.17}
Let $ X : \Om \to \R/\Q $ have the absolutely continuous distribution
(recall \ref{2.32}). Then $ \Pr{ X \in A_1, X \in A_2 } = \Pr{ X \in A_1
} \Pr{ X \in A_2 } $ for all $ A_1,A_2 \in \Si $, but $ X $ is not
independent of itself.
\end{counterexample}

\begin{proof}
The equality $ \Pr{ X \in A_1, X \in A_2 } = \Pr{ X \in A_1 } \Pr{ X \in
A_2 } $ is easy to check in each of the four possible cases ($ 0\cdot0
$, $ 0\cdot1 $, $ 1\cdot0 $, $ 1\cdot1 $) taking into account that $ X $
has a \OI\ distribution. (See also \cite[4.6]{Ke00}.) Assume that $ X =
\Q + Y_1 $ and $ X = \Q + Y_2 $ where $ Y_1,Y_2 : \Om \to \R $ are
independent. Then $ \Pr{ Y_1 - Y_2 \in \Q } = 1 $, therefore $ \Pr{ Y_1
- Y_2 = q } > 0 $ for some $ q \in \Q $, and $ 0 < \cP{ Y_1 - Y_2 = q }{
Y_2 = y_2 } = \Pr{ Y_1 = y_2+q } $ for some $ y_2 $; that is, $ Y_1 $
has an atom at $ y_1 = y_2 + q $. On the other hand, $ X_1 = \Q + Y $
for some absolutely continuous $ Y $. We have $ \Pr{ Y_1 - Y \in \Q } =
1 $, therefore $ \Pr{ Y \in \Q + y_1 } > 0 $ in contradiction to the
absolute continuity of $ Y $.
\end{proof}

\section[]{\raggedright Some generalizations and final remarks}
\label{sect9}\begin{remark}
Theorem \ref{6n1} remains true if `the independence condition' is
replaced with `the quasi-independence condition' introduced as follows.

\begin{sloppypar}
First, in Definition \ref{4n1}, instead of item (c) ``the random
sequence $ (Y_2,Y_4,Y_6,\dots) $ is independent of the random sequence $
(Y_1,Y_3,Y_5,\dots) $'' we write ``the \sif s $ \si(Y_1,Y_3,Y_5,\dots)
$, $ \si(Y_2,Y_4,Y_6,\dots) $ are a nonsingular pair (on the whole
probability space)'', in which case we call the two fragments
\emph{quasi-independent}. Recall that `nonsingular pair' means
(according to Definition \ref{5n5}) absolute continuity of the joint
distribution w.r.t.\ the product of its two marginals.
\end{sloppypar}

Second, in Definition \ref{4n2}, instead of ``the $ n $ fragments \dots\
are independent'' we write ``the $ n $ fragments \dots\ are
quasi-independent'', thus defining \emph{the quasi-independence
condition.}

The only change in the proof of Proposition \ref{5n2} is, deletion of
the last four words `since they are independent'. The proof of Lemma
\ref{5n10} still holds. Now, Propositions \ref{5n1}, \ref{5n9} and
Theorem \ref{6n1} use these modified \ref{5n2} and \ref{5n10} as
before.

We see that quasi-independence is no less restrictive than independence,
for random dense countable sets. For random \emph{closed} sets the
situation is completely different; independence characterizes Poisson
processes, while quasi-independence leaves a great freedom
\cite[Sect.~6]{Ts02}.
\end{remark}

\begin{remark}\label{9n2}
Returning to independence, we compare our approach with that of Kingman
\cite[Sect.~2]{Ki}. A Poisson process is defined there as a family $ \(
X(\om) \)_{\om\in\Om} $ of finite or countable subsets $ X(\om) $ of a
given measurable space $ T $, satisfying three conditions: measurability,
independence and distribution. The \emph{measurability} condition: for
every measurable $ B \subset T $ the number $ \xi_B(\om) \in
\{0,1,2,\dots\} \cup \{\infty\} $ of points in $ B \cap X(\om) $ is
measurable in $ \om $. The \emph{independence} condition: for any
disjoint measurable $ B_1,\dots,B_n \subset T $ the random variables $
\xi_{B_1}, \dots, \xi_{B_n} $ are independent. The \emph{distribution}
condition: each random variable $ \xi_B $ has the Poisson distribution
with some parameter $ \mu(B) \in [0,\infty] $; here $ \mu(B)=0 $ means
that $ \xi_B = 0 $ a.s., while $ \mu(B)=\infty $ means that $ \xi_B =
\infty $ a.s.

Clearly, $ \mu $ is a measure, positive, maybe infinite and even not
\finite{\si}, and nonatomic in the sense that $ \mu(\{t\}) = 0 $ for all
$ t \in T $. It appears \cite[Sect.~2.5]{Ki} that such a Poisson process
exists for every $ \mu $ of the form $ \mu_1 + \mu_2 + \dots $ where $
\mu_n $ are nonatomic \emph{finite} positive measures on $ T $.

For example, we may take $ T = (0,1) $ and $ \mu_1 = \mu_2 = \dots $ be
Lebesgue measure on $ (0,1) $. Then Kingman's construction gives just
the object $ X $ that we call an unordered infinite sample. However,
consider $ Y(\om) = ( \Q + y(\om) ) \cap (0,1) $; here $ \Q \subset \R $
is the set of all rational numbers, $ \Q + y(\om) $ its shift by $
y(\om) $, and $ y : \Om \to \R $ a random variable with an absolutely
continuous distribution (as in \ref{2.32}, \ref{2.17}). Is $ Y $ also a
Poisson process?

In the framework of Kingman, $ Y $ is a Poisson process, and moreover, $
X $ and $ Y $ are treated as identically distributed, just because
random variables $ \xi_B $ do not feel any difference between $ X $ and
$ Y $. However, distances between points are irrational for $ X $ but
rational for $ Y $; a clear-cut distinction!

In our framework (recall \ref{2m1} and \ref{2.3}) $ X $ and $ Y $ are
not identically distributed, and $ Y $ should not be called a Poisson
process, since it violates the independence condition \ref{4n2}
(similarly to \ref{2.17}).

See also \cite[Sect.~2.2]{Ki}: ``It might be objected that $ \Pi_1 $ and
$ \Pi_2 $ are not `really' independent, and only appear to be so because
we choose to describe them in terms of their count processes.''
\end{remark}

\begin{counterexample}
Kingman's measurability condition (mentioned in \ref{9n2}) does not
imply the measurability condition \eqref{2n0}. There exists a family $
\( X(\om) \)_{\om\in\Om} $ of countable sets $ X(\om) \subset \R $ not
of the form \eqref{2n0} but such that all $ \xi_B $ are measurable. Here
$ \Om $ is a standard probability space.

We take $ \Om = (0,1) $ (with Lebesgue measure), choose an irrational
number $ a \in \R $ and a set $ A \subset (0,1) $ of interior measure $
0 $ but outer measure $ 1 $, and define $ X $ by
\[
X(\om) = \begin{cases}
 \Q + \om &\text{for $ \om \in A $},\\
 \Q + \om + a &\text{for $ \om \in (0,1) \setminus A $};
\end{cases}
\]
here $ \Q $ is the set of all rational numbers, and $ \Q + \om $ its
shift by $ \om $.

If a Borel set $ B \subset \R $ is negligible (that is, of Lebesgue
measure $ 0 $) then $ \{ \om : X(\om) \cap B \ne \emptyset \} $ is
negligible (since it has negligible intersection with $ A $ and also
with $ (0,1) \setminus A $); thus $ \xi_B = 0 $ a.s. Otherwise, if $ B $
is of positive measure, then $ \{ \om : X(\om) \cap B = \emptyset \} $
is negligible (for the same reason); it follows that $ \xi_B = \infty $
a.s.

Assume that a function $ X_1 : (0,1) \to \R $ is such that $ X_1(\om)
\in X(\om) $ for almost all $ \om $. Then $ X_1(\om) - \om $ is rational
for almost all $ \om \in A $ but irrational for almost all $ \om \in
(0,1) \setminus A $. Thus $ X_1 $ cannot be measurable, which shows that
$ X $ is not of the form \eqref{2n0}.
\end{counterexample}

Kingman's independence condition is too demanding for the
percolation-related models mentioned in Introduction; independent
fragments of their `full scaling limit' are well-defined over disjoint
good domains, not just Borel sets. An appropriate independence condition
is given below.

Random finite or countable subsets of a standard Borel space $ T $ are
random elements (see \ref{8n1}) of the borelogical space $ \FCS(T) $
defined similarly to \ref{7n8}. They may be written as $ X = \{
X_1,\dots,X_N \} $ where $ X_1,X_2,\dots : \Om \to T $ and $ N : \Om \to
\{0,1,2,\dots\} \cup \{\infty\} $ are random variables. Here $ X(\om) =
\{ X_1(\om), \dots, X_{N(\om)}(\om) \} $ when $ N(\om) < \infty $ and $
X(\om) = \{ X_1(\om), X_2(\om), \dots \} $ when $ N(\om) = \infty $. Of
course, $ N(\om) = 0 $ means $ X(\om) = \emptyset $.

Given a Borel set $ B \subset T $, the fragment $ \om \mapsto B \cap
X(\om) $ of $ X $ is again a random finite or countable
set. Independence of two or more such fragments is understood according
to \ref{2.16}. It is less clear how to generalize \ref{4n2} since it is
based on intervals.

\begin{definition}\label{9n3}
Let $ T $ be a standard Borel space, $ X $ a random finite or countable
subset of $ T $, and $ \A $ an algebra of Borel subsets of $ T $. We say
that $ X $ satisfies \emph{the independence condition} on $ \A $, if for
every $ n = 2,3,\dots $ and every $ n $ disjoint sets $ B_1, \dots, B_n
\in \A $, the $ n $ fragments $ B_1 \cap X, \dots, B_n \cap X $ are
independent.
\end{definition}

By the way, a \emph{countable} algebra $ \A $ generates the Borel \sif\
if and only if it separates points (see \cite[(14.16)]{Ke95}).

\begin{lemma}\label{9n4}
If $ X $ satisfies the independence condition on some algebra $ \A $
that generates the Borel \sif\ of $ T $, and $ \Pr{ t \in X } = 0 $ for
all $ t \in T $, then $ X $ is a Poisson process in the sense of Kingman
(see \ref{9n2}).
\end{lemma}

\begin{proof}
The three conditions mentioned in \ref{9n2} must be verified. The
measurability condition holds evidently. The independence condition
evidently holds for $ B_1,\dots,B_n \in \A $; after some preparations it
will be generalized to all measurable $ B_1,\dots,B_n $.

There exists a nonatomic finite positive measure $ \nu $ on $ T $ such
that
\[
\nu(B) = 0 \quad \text{if and only if} \quad \Pr{ B \cap X \ne \emptyset
} = 0
\]
for all Borel $ B \subset T $. For example, one may represent $ X $ as $
\{ X_1,\dots,X_N \} $ and take $ \nu(B) = \sum_n 2^{-n} \Pr{ N \ge n,
X_n \in B } $.

For every Borel $ B \subset T $ there exist $ A_1,A_2,\dots \in
\A $ such that $ \nu ( B \bigtriangleup A_n ) \to 0 $ and moreover, $
\sum_n \nu ( B \bigtriangleup A_n ) < \infty $. (Of course, $ B
\bigtriangleup A_n = ( B \setminus A_n ) \cup ( A_n \setminus B ) $.)
We introduce sets $ C_n = \bigcap_{k\ge n} A_k $, $ C = \bigcup_n C_n $
and $ D_n = \bigcap_{k\ge n} ( T \setminus A_k ) $, $ D = \bigcup_n D_n
$, then $ C_n \uparrow C $, $ D_n \uparrow D $, $ C \cap D = \emptyset $
and $ \nu ( B \setminus C ) = 0 $, $ \nu ( ( T \setminus B ) \setminus D
) = 0 $ by the first Borel-Cantelli lemma. Thus, $ \nu ( B
\bigtriangleup C ) = 0 $ and $ \nu ( ( T \setminus B ) \bigtriangleup D
) = 0 $.

For each $ n $ we have $ C_n \subset A_n $ and $ D_n \subset T \setminus
A_n $; the independence condition on $ \A $ implies that $ \xi_{C_n},
\xi_{D_n} $ are independent. However, $ \xi_{C_n} \uparrow \xi_C = \xi_B
$ and $ \xi_{D_n} \uparrow \xi_D = \xi_{T \setminus B} $; we conclude
that $ \xi_B $ and $ \xi_{T \setminus B} $ are independent for every
Borel set $ B \subset T $. Similarly, $ \xi_{B_1}, \dots, \xi_{B_n} $
are independent whenever Borel sets $ B_1,\dots,B_n $ are disjoint (this
generalization is left to the reader); the independence condition is
verified.

Without loss of generality we assume that $ T = (0,1) $ and $ \nu $ is
absolutely continuous w.r.t.\ Lebesgue measure (by the isomorphism
theorem for measure spaces, see \cite[(17.41)]{Ke95} or
\cite[3.4.23]{Sri}). We cannot apply Prop.~\ref{4n4} `as is', since the
algebra $ \A $ need not contain intervals (essential in
\ref{4n2}). However, the independence condition \ref{4n2} is used in the
proof of \ref{4n4} only once; it ensures independence of the random
variables $ \xi(x,y) $ for disjoint intervals $ (x,y) $. In our case,
independence between these $ \xi(x,y) $ is ensured by the independence
condition formulated in \ref{9n2} and verified above. Thus, the
conclusion of Prop.~\ref{4n4} holds; the distribution condition is
verified.
\end{proof}

\begin{remark}\label{9n5}
\begin{sloppypar}
Let $ X $ be as in Lemma \ref{9n4}, then $ \mu $ defined by $ \xi_B \sim
\Poisson(\mu(B)) $ is a measure of the form $ \mu_1 + \mu_2 + \dots $
where $ \mu_n $ are nonatomic \emph{finite} positive measures on $ T $.
Indeed, representing $ X $ as $ \{ X_1,\dots,X_N \} $ we have $ \mu(B) =
\sum_n \Pr{ N \ge n, X_n \in B } $.
\end{sloppypar}

On the other hand, every such measure $ \mu $ corresponds to some $ X $,
which follows from Kingman's construction (mentioned in \ref{9n2}).
\end{remark}

\begin{remark}\label{9n6}
The independence condition \ref{4n2} may be treated as a special case of
the independence condition \ref{9n3}; namely, $ T = (0,1) $ and $ \A $
consists of finite unions of intervals. Similarly we may take $ T = \R^2
$ and $ \A $ consisting of finite unions of rectangles. Using \ref{9n4}
we may generalize \ref{4n4} to the two-dimensional case. This way, the
main results (Theorems \ref{6n1}, \ref{6n8}) may be generalized to
dimension $ 2 $ (and higher).
\end{remark}

\begin{counterexample}\label{9n75}
It may happen that the distribution of the sequence $
(X_2,X_1,X_4,X_3,X_6,X_5,\dots) $ satisfies the conditions of Main lemma
\ref{3.2}, but the distribution of the sequence $ (X_1,X_2,X_3,\dots) $
does not.

We construct random variables $ X_k $ via their binary digits $
\be_{k,l} : \Om \to \{0,1\} $,
\[
X_k = (0.\be_{k,1}\be_{k,2}\dots)_2 = \sum_l 2^{-l} \be_{k,l} \, .
\]
Each $ \be_{k,l} $ takes on the two values $ 0,1 $ with the
probabilities $ \frac12, \frac12 $ and they all are independent except
for the following restriction:
\[
\be_{k+1,l} = \be_{k,k+l-1} \quad \text{for } l = 1,\dots,k
\]
and $ k = 1,2,\dots $

The conditional distribution of $ X_{n+1} $ given $ X_1,\dots,X_n $ is
the uniform distribution on an interval of length $ 2^{-n} $. Thus,
Condition \ref{3.2}(a) is satisfied, but \ref{3.2}(b) is violated;
moreover, the series of \ref{3.2}(b) converges almost everywhere.

Condition \ref{3.2}(a), being permutation-invariant, is still satisfied
by the distribution of $ (X_2,X_1,X_4,X_3,X_6,X_5,\dots) $. Condition
\ref{3.2}(b) is also satisfied, since $ X_{2n} $ is independent of $
X_1,\dots,X_{2n-2} $, which makes every second term of the series equal
to $ 1 $ almost everywhere.
\end{counterexample}

\begin{remark}\label{9n8}
Let $ \mu_1,\mu_2 $ be two different probability measures on $ (0,1) $,
equivalent (that is, mutually absolutely continuous) to Lebesgue
measure. Consider the mixture $ \nu = \frac12 ( \mu_1^\infty +
\mu_2^\infty ) $ of the corresponding product measures $ \mu_1^\infty,
\mu_2^\infty $ on $ (0,1)^\infty $. The measure $ \nu $ is invariant
under the group $ S^\infty $ of all permutations, that is, invertible
maps $ s : \{ 1,2,\dots \} \to \{ 1,2,\dots \} $ ($ S^\infty $ acts on $
(0,1)^\infty $ by $ sx = ( x_{s(1)}, x_{s(2)}, \dots ) $ for $
x=(x_1,x_2,\dots) $). Consider also the countable subgroup $ S_\infty
\subset S^\infty $ consisting of $ s $ such that the set $ \{ n : s(n)
\ne n \} $ is finite.

There exists an \invariant{S_\infty} Borel set $ B \subset (0,1)^\infty
$ such that $ \mu_1^\infty (B) = 1 $ but $ \mu_2^\infty (B) = 0 $. For
example we may choose a Borel $ A \subset (0,1) $ such that $ \mu_1(A) <
\mu_2(A) $ and take
\[
B = \Big\{ (x_1,x_2,\dots) : \limsup_{n\to\infty} \frac1n \(
\One_A(x_1) + \dots + \One_A(x_n) \) < \mu_2(A) \Big\} \, .
\]
The set $ B $ is not \invariant{S^\infty}, however, it is
\invariant{S^\infty} $ \modO $, that is, $ \nu ( B \bigtriangleup sB ) =
0 $ for each $ s \in S^\infty $.

In contrast, $ \mu_1^\infty (B) = \mu_2^\infty (B) $ for every
\invariant{S^\infty} Borel set $ B \subset (0,1)^\infty $. \emph{Proof.}
Both $ \mu_1^\infty $ and $ \mu_2^\infty $ satisfy the conditions of
Main lemma \ref{3.2}, therefore they lead to the same distribution on $
\CS(0,1) $ (recall \ref{2n2}).
\end{remark}

\begin{remark}
Finite or countable sets may be treated as equivalence classes of
discrete probability measures (equivalence being mutual absolute
continuity), see \ref{7n8}. Equivalence classes of nonatomic
\emph{singular} measures are another borelogical space. Random elements
of this space may be subjected to conditions of independence and
stationarity. An interesting example associated with Brownian motion is
discussed in \cite[Sect.~2f]{Ts04} in connection with a nonclassical
noise (Warren's noise of stickiness).
\end{remark}

\begin{remark}
Independence and stationarity of Brownian local minimizers result from
(a) independence and stationarity of Brownian increments (the white
noise) and (b) factorizability and stationarity of the map from Brownian
increments to Brownian local minimizers. In terms of
\cite[Sect.~2e]{Ts04} this map is an example of a stationary local
random dense countable set over the white noise. Brownian local
maximizers are another example. Their union, Brownian local extrema, are
the third example. (Several types of special points on the Brownian path
are examined, see \cite{Pe}, but they are uncountable sets.) The
question \cite[2e3]{Ts04}, whether or not these three examples exhaust
\emph{all} stationary local random dense countable sets over the white
noise, is still open!
\end{remark}

\bigskip
\filbreak
{
\small
\begin{sc}
\parindent=0pt\baselineskip=12pt
\parbox{4in}{
Boris Tsirelson\\
School of Mathematics\\
Tel Aviv University\\
Tel Aviv 69978, Israel
\smallskip
\par\quad\href{mailto:tsirel@post.tau.ac.il}{\tt
 mailto:tsirel@post.tau.ac.il}
\par\quad\href{http://www.tau.ac.il/~tsirel/}{\tt
 http://www.tau.ac.il/\textasciitilde tsirel/}
}

\end{sc}
}
\filbreak

\end{document}